\documentclass[11pt]{amsart}

\usepackage{amssymb,amsmath,amsthm}
\usepackage{type1cm}
\usepackage{xcolor}
\usepackage{graphicx}
\usepackage{multicol,booktabs}
\usepackage[margin=1in]{geometry}
\usepackage[font=small,margin=0pt,belowskip=0pt,aboveskip=5pt,tableposition=top]{caption}
\usepackage{algorithm,algorithmic}
\usepackage{nicefrac}

\usepackage[backref=true,sorting=none,citestyle=numeric-comp]{biblatex}
\usepackage[inline]{enumitem}
\usepackage{siunitx}

\definecolor{red}{RGB}{200,16,46}

\usepackage[colorlinks=true,linkcolor=red,citecolor=red]{hyperref}

\newcommand{\figref}[1]{Figure~\ref{#1}}
\newcommand{\tabref}[1]{Table~\ref{#1}}
\newcommand{\secref}[1]{\S\ref{#1}}
\newcommand{\algref}[1]{Algorithm~\ref{#1}}

\newcommand\fs[1]{\ensuremath{\mathfrak{#1}}}
\newcommand\op[1]{\ensuremath{\mathcal{#1}}}
\newcommand\di[1]{\ensuremath{\mathbf{#1}}}
\newcommand\igrad{\ensuremath{\nabla}}
\newcommand{\idiv}{\ensuremath{\nabla\cdot}}
\newcommand{\ilap}{\rotatebox[origin=c]{180}{$\nabla$}}

\newcommand{\defeq}{\ensuremath{\mathrel{\mathop:}=}}

\DeclareMathOperator*\minopt{\ensuremath{minimize}}
\newcommand{\reg}{\ensuremath{\text{reg}}}

\newcommand{\dist}{\ensuremath{\text{dist}}}
\newcommand{\st}{\ensuremath{\text{subject to}}}
\renewcommand{\d}{\ensuremath{\text{d}}}
\renewcommand{\vec}[1]{\ensuremath{\boldsymbol{#1}}}

\newcommand\dual{\ensuremath{\ast}}

\newcommand\fnum[1]{\num[scientific-notation=fixed,round-precision=3,round-mode=places]{#1}}
\newcommand\snum[1]{\num[round-precision=3,round-mode=places,tight-spacing=true,retain-zero-exponent=true,retain-explicit-plus=true,output-exponent-marker=\text{e}]{#1}}
\newcommand\inum[1]{\num[scientific-notation=fixed,round-precision=0,round-mode=places,group-separator = {,}]{#1}}
\newcommand\scinum[1]{\num[round-precision=0,round-mode=places,tight-spacing=true,retain-zero-exponent=true,retain-explicit-plus=true,output-exponent-marker=\text{e}]{#1}}

\title[CLAIRE: GPU-Accelerated Algorithms for Diffeomorphic Image Registration]{CLAIRE: Scalable GPU-Accelerated Algorithms for Diffeomorphic Image Registration in 3D}

\author{Andreas Mang}
\address{Department of Mathematics, University of Houston, 3551 Cullen Blvd, Houston, TX 77204-3008}
\email{andreas@math.uh.edu}

\begin{document}

\maketitle

\begin{abstract}
We present our work on scalable, GPU-accelerated algorithms for diffeomorphic image registration. The associated software package is termed CLAIRE. Image registration is a non-linear inverse problem. It is about computing a spatial mapping from one image of the same object or scene to another. In diffeomorphic image registration, the set of admissible spatial transformations is restricted to maps that are smooth, one-to-one, and have a smooth inverse. We formulate diffeomorphic image registration as a variational problem governed by transport equations. We use an inexact, globalized (Gauss--)Newton--Krylov method for numerical optimization. We consider semi-Lagrangian methods for numerical time integration. Our solver features mixed-precision, hardware-accelerated computational kernels for optimal computational throughput. We use the message-passing interface for distributed-memory parallelism and deploy our code on modern high-performance computing architectures. Our solver allows us to solve clinically relevant problems in under four seconds on a single GPU. It can also be applied to large-scale 3D imaging applications with data that is discretized on meshes with billions of voxels. We demonstrate that our numerical framework yields high-fidelity results in only a few seconds, even if we search for an optimal regularization parameter.
\end{abstract}

\section{Introduction}
\label{s:intro}

\begin{figure}
\centering
\includegraphics[width=\textwidth]{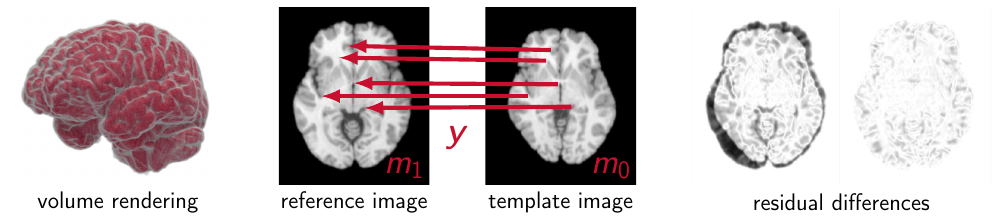}
\caption{Image registration problem. On the left, we show a volume rendering of a 3D brain MRI. The figures in the middle show an axial slice of two MRI brain scans of different individuals. In image registration, we seek a map $\vec{y} \in \fs{Y}_{\text{ad}} \subset \{\vec{\phi} : \mathbb{R}^d \to \mathbb{R}^d\}$, $d\in\{2,3\}$, that establishes a \emph{plausible} spatial correspondence between these to images. In this work, we restrict the set of admissible spatial transformations $\fs{Y}_{\text{ad}}$ to $\mathbb{R}^d$-diffeomorphisms. On the right, we show the residual differences between the axial slices of the images shown in the middle before (left) and after (right) diffeomorphic (deformable) registration. Here, white represents small residual differences, and black indicates large residual differences. We note that the registration of two brains from different individuals is a common application for diffeomorphic image registration in computational anatomy. However, strictly speaking this example is in violation with our underlying assumptions; we do not register the ``same object''---we register images of brains of different individuals (in an attempt to study anatomical variability).\label{f:imagereg}}
\end{figure}

In the present work, we discuss scalable, hardware-accelerated algorithms for diffeomorphic image registration. We review our past contributions and showcase results for a software framework termed CLAIRE~\cite{Mang:2019a,claireweb,Brunn:2021b}. Image registration is an ill-posed inverse problem~\cite{Fischer:2008a}. It is a key methodology in medical image analysis. The inputs are two (or more, noisy) images $m_i \in \fs{I}$, $i = 0,1$, $\fs{I} \subset \{u : \Omega \to \mathbb{R}\}$, of the same object or scene, compactly supported on some domain $\Omega \subseteq \mathbb{R}^d$, where $d \in \{2,3\}$. In image registration, we seek a \emph{plausible} spatial transformation $\vec{y} \in \fs{Y}_{\text{ad}}$, $\fs{Y}_{\text{ad}} \subset \{\vec{\phi} : \mathbb{R}^d \to \mathbb{R}^d\}$, that maps points in the so-called \emph{template} or \emph{source image} $m_0$ to its corresponding points in the so-called \emph{reference} or \emph{target image} $m_1$~\cite{Modersitzki:2004a,Modersitzki:2009a,Fischer:2008a}. The notion of the plausibility of the map $\vec{y}$ depends on the particular application. In the present work, we restrict the set $\fs{Y}_{\text{ad}}$ of admissible maps $\vec{y}$ to $\mathbb{R}^d$-diffeomorphisms~\cite{Younes:2019a}. That is, $\fs{Y}_{\text{ad}} \subseteq \operatorname{diff}(\mathbb{R}^d)$, where $\operatorname{diff}(\mathbb{R}^d)$ is the  set of $\mathbb{R}^d$-diffeomorphisms, i.e., smooth maps from $\mathbb{R}^d$ to $\mathbb{R}^d$ that are one-to-one and onto, with a smooth inverse. The set $\operatorname{diff}(\mathbb{R}^d)$ is closed under composition and taking the inverse; it forms a group. In this framework, deforming the template image $m_0$ corresponds to a change of coordinates $m_0 \circ \vec{y}^{-1}$; the image intensity in the transformed image $m_0 \circ \vec{y}^{-1}$ at coordinate $\vec{y}(\vec{x})\in\mathbb{R}^d$ is identical to the value at the location $\vec{x}\in\mathbb{R}^d$ in the original image. Using this notation, the diffeomorphic image registration problem can be formulated as the problem of finding $\vec{y} \in \operatorname{diff}(\mathbb{R}^d)$ such that $m_0 \circ \vec{y}^{-1} = m_1$. We illustrate this in~\figref{f:imagereg}. We summarize the main notation and acronyms in \tabref{t:notation}.

\begin{table}
\caption{Notation, symbols, and acronyms.\label{t:notation}}
\centering
\begin{tabular}{ll}\toprule
\bf Symbol/Acronym & \bf Meaning \\
\midrule
$d \in \mathbb{N}$ & dimensionality of the ambient space\\
$\Omega \subset \mathbb{R}^d$ & spatial domain \\
$m_0 : \bar{\Omega} \to \mathbb{R}$ & template image \\
$m_1 : \bar{\Omega} \to \mathbb{R}$ & reference image \\
$m : [0,1] \times \bar{\Omega} \to \mathbb{R}$ & state variable (transported image intensities) \\
$\vec{v} : \bar{\Omega} \to \mathbb{R}^d$ & control variable (stationary velocity field)\\
$\lambda : [0,1] \times \bar{\Omega} \to \mathbb{R}$ & dual variable \\
$\dist : \fs{I} \times \fs{I} \to \mathbb{R}$ & distance functional \\
$\reg : \fs{V} \to \mathbb{R}$ & regularization functional \\
$\op{L} : \fs{V} \to \fs{V}^\dual$ & regularization operator \\
$\alpha \in \mathbb{R}$ & regularization parameter \\
$\fs{I}$ & orbit \\
$\fs{G}$ & group of diffeomorphisms \\
$\operatorname{diff}(\mathbb{R}^d)$ & set of $\mathbb{R}^d$-diffeomorphisms \\
\midrule
FFT   & Fast Fourier Transform\\
GMRES & Generalized Minimal RESidual (method)\\
GPU   & Graphics Processing Unit\\
HPC   & High Performance Computing\\
KKT   & Karush--Kuhn--Tucker (conditions)\\
LDDMM & Large Deformation Diffeomorphic Metric Mapping\\
MRI   & Magnetic Resonance Imaging\\
MPI   & Message Passing Interface\\
ODE   & Ordinary Differential Equation\\
PCG   & Preconditioned Conjugate Gradient (method)\\
PDE   & partial differential equation\\
RK2   & second-order Runge--Kutta (method)\\
SL    & Semi-Lagrangian (method)\\
\bottomrule
\end{tabular}
\end{table}

\subsection{Outline of the Method}

The approach considered in the present work is related to a mathematical framework referred to as LDDMM~\cite{Younes:2019a,Beg:2005a,Trouve:1995a,Trouve:1998a,Dupuis:1998a}. We consider PDE-constrained optimization problems~\cite{Borzi:2012a,Hinze:2009a,Antil:2018a,Mang:2018a} governed by transport equations for diffeomorphic image registration. The transport map is parameterized by a smooth space-time field $\vec{v} \in \fs{V}$, $\fs{V} \defeq L^q([0,1],\fs{H})$, $q\in \mathbb{N}$, where $\fs{H}$ is a Sobolev space of suitable regularity, i.e., $\fs{H} = W^{p,s}(\Omega,\mathbb{R}^d)$, $p\in\mathbb{N}$, $s\in\mathbb{N}$. Our problem formulation is of the form
\begin{subequations}
\label{e:varopt}
\begin{align}
\minopt_{m \in \fs{M}_{\text{ad}}, \, \vec{v} \in \fs{V}_{\text{ad}}}\quad
&\dist(m(1), m_1) + \reg(\vec{v})
\label{e:varopt:objective}
\\
\st\quad & c(m,\vec{v}) = 0.
\label{e:varopt:constraint}
\end{align}
\end{subequations}

\noindent Here, $c : \fs{M} \times \fs{V} \to \fs{Q}$ represents a PDE constraint. It is of the general form $c(m,\vec{v}) = \op{A}(m,\vec{v}) - q$. The \emph{parameter-to-observation map} $f : \fs{V} \to \fs{I}$ (i.e., the solution operator for the constraint) is formally given by $f(\vec{v}) = \op{Q} \op{A}^{-1}(m,\vec{v})q$. Here, $\op{Q}$ denotes the observation operator, i.e., a mapping that takes the output of $\op{A}^{-1}$ and maps it to ``locations'' at which data is available. The functional $\dist : \fs{I} \times \fs{I} \to \mathbb{R}$ measures the discrepancy between the deformed template image $m(1)= f(\vec{v})$ and the reference image $m_1$. The functional $\reg : \fs{V} \to \mathbb{R}$ is a regularization functional. We specify the precise choices in greater detail below.

We use the method of Lagrange multipliers to solve~\eqref{e:varopt}. We consider an \emph{optimize-then-discretize} approach. We use a globalized, inexact reduced space \mbox{(Gauss--)}Newton--Krylov method for numerical optimization. We solve the PDEs that appear in the optimality conditions based on a SL method. The main computational kernels of our algorithm are interpolation and numerical differentiation. For interpolation, we use a Lagrange polynomial. For numerical differentiation, we consider a mixture of high-order finite difference operators and a pseudo-spectral method. We use MPI for distributed-memory parallelism and deploy our code on dedicated GPU architectures.

\subsection{Related Work}

We consider a PDE-constrained optimization problem for velocity-based diffeomorphic image registration. We refer to~\cite{Biegler:2003a,Borzi:2012a,Gunzburger:2003a,Hinze:2009a,Antil:2018a,Lions:1971a} for insights into theory and algorithmic developments related to PDE-constrained optimization. Additional information about image registration and related work can be found in~\cite{Modersitzki:2004a,Fischer:2008a,Modersitzki:2009a,Hajnal:2001a,Sotiras:2013a,Younes:2019a}. As we mentioned above, we restrict ourselves to diffeomorphic image registration. An intuitive approach to safeguard against non-diffeomorphic maps $\vec{y}$ is to add hard and/or soft constraints to the variational problem~\cite{Burger:2013a,Haber:2007a,Rohlfing:2003a,Sdika:2008a}. An alternative strategy is to introduce a pseudo-time variable $t$ and invert for a smooth velocity field $\vec{v}$ that parameterizes $\vec{y}$~\cite{Younes:2019a,Younes:2020a,Beg:2005a,Dupuis:1998a,Trouve:1998a,Miller:2001a,Vercauteren:2009a,Christensen:1996a,Younes:2007a}; our approach falls into this category. In~\cite{Beg:2005a,Dupuis:1998a,Trouve:1995a,Trouve:1998a,Miller:2001a,Younes:2007a}, the flow of the sought after diffeomorphism $\vec{y}$ is modelled as the solution of the ODE $\partial_t \vec{\phi} = \vec{v} \circ \vec{\phi}$ for $t \in (0,1]$ with initial condition $\vec{\phi} = \operatorname{id}_{\mathbb{R}^d}$ at time $t = 0$, where $\vec{v}$ is a smooth, time-dependent vector field from $\mathbb{R}^d$ to $\mathbb{R}^d$ and $\operatorname{id}_{\mathbb{R}^d} : \mathbb{R}^d \to \mathbb{R}^d$, $\operatorname{id}_{\mathbb{R}^d}(\vec{x}) = \vec{x}$, is the identity transformation in $\mathbb{R}^d$. This ODE enters the variational problem as a constraint; we arrive at a non-linear optimal control problem with state variable $\vec{\phi}$ and control $\vec{v}$. The sought-after diffeomorphism $\vec{y}$ that maps one image to another corresponds to the end point of the flow $\vec{\phi}$, i.e., $\vec{y} = \vec{\phi}(t=1)$. This approach is commonly referred to as LDDMM~\cite{Beg:2005a}. We describe it in greater detail in the main part of this manuscript. In our formulation, the diffeomorphism $\vec{y}$ does no longer appear; we do \emph{not} model the deformed template image as the application of $\vec{y}$ to $m_0$. Instead, we transport the intensities of the template image $m_0$ given some candidate $\vec{v}$ based on a hyperbolic transport equation~\cite{Mang:2015a,Hart:2009a,Borzi:2002a}. Unlike most existing approaches, our framework features explicit control on volume change introduced by the mapping by controlling the divergence of $\vec{v}$. This formulation was originally proposed in~\cite{Mang:2016a}; a similar approach is described in~\cite{Borzi:2002a}. Works of other groups that consider divergence-free velocities $\vec{v}$ in similar contexts have been described in~\cite{Chen:2011a,Hinkle:2009a,Mansi:2011a,Ruhnau:2007a,Saddi:2008a}.

Our formulation has been introduced in~\cite{Mang:2015a,Mang:2016a}. The work most closely related to ours in terms of the problem formulation is~\cite{Borzi:2002a,Borzi:2002b,Chen:2012a,Hart:2009a,Lee:2010a,Lee:2011a,Vialard:2012a,Herzog:2019a}. Related formulations for optimal mass transport are discussed in~\cite{Benzi:2011a,Haber:2015a,Rehman:2009a,Mang:2017a,Herzog:2019a}. In contrast to optimal mass transport, our formulation keeps the transported quantities constant along the characteristics, i.e., mass is \emph{not} preserved. Our formulation is related to traditional optical flow formulations~\cite{Horn:1981a,Kalmoun:2011a,Ruhnau:2007a}. The main difference is that the transport equation for the image intensities of $m_0$ enters our formulation as a hard constraint. PDE-constrained formulations for optical flow that are equivalent to our formulation are described in~\cite{Andreev:2015a,Barbu:2016a,Borzi:2002a,Chen:2011a}.

Among the most popular packages for diffeomorphic registration are {\tt Demons}~\cite{Vercauteren:2008a,Vercauteren:2009a}, {\tt ANTs}~\cite{Avants:2011a,Avants:2008a,}, {\tt Deformetrica}~\cite{Bone:2018a,Fishbaugh:2017a}, or {\tt DARTEL}~\cite{Ashburner:2007a}. There are only few works on effective numerical methods for velocity-based diffeomorphic image registration, and even fewer on scalable algorithms. Works of other groups on numerical algorithms for the solution of diffeomorphic  association problems (for images as well as surface representations) are, e.g., described in~\cite{Beg:2005a,Cao:2005a,Hsieh:2022a,Polzin:2016a,Polzin:2020a,Vialard:2012a,Niethammer:2009a,Hart:2009a,Arguillere:2016a}. The majority of existing works consider an \emph{optimize-then-discretize} approach for solving the variational problem~\cite{Beg:2005a,Arguillere:2016a,Vialard:2012a,Zhang:2015b,Ashburner:2011a,Hart:2009a,Miller:2006a}; \emph{discretize-then-optimize} approaches for related problem formulations can be found in~\cite{Polzin:2020a,Azencott:2010a,Mang:2017a,Zhang:2021a,Mang:2023a}. In the work discussed in this exposition, we also consider an \emph{optimize-then-discretize} approach~\cite{Mang:2015a,Mang:2016a,Mang:2016b}; an implementation for a \emph{discretize-then-optimize} approach for problem formulations similar to the one considered here can be found in~\cite{Mang:2017a}.

Despite the fact that first-order methods for optimization have poor convergence rates for nonlinear, ill-posed inverse problems, most work on algorithms for formulations similar to ours, with the exception of ours~\cite{Mang:2015a,Mang:2016a,Mang:2018a,Mang:2016b,Mang:2017a,Mang:2017b} and~\cite{Ashburner:2011a,Benzi:2011a,Hernandez:2014a,Simoncini:2012a,Vercauteren:2009a,Herzog:2019a,Polzin:2016a}, use first order gradient descent-type approaches. Work on operator-splitting algorithms for LDDMM (and related problems) can be found in~\cite{Thorley:2021a,Lee:2016a,Zhang:2021a,Mang:2023a}. Other recent works that do not explicitly derive optimality conditions based on variational principles but rely on automatic differentiation can be found in~\cite{Hsieh:2021a,Franccois:2021a,Bone:2018a,Hartman:2023a,Bone:2020a}. Lastly, we note that the success of machine learning in various scientific disciplines has led to several recent works that attempt to solve the inverse problem of diffeomorphic registration based on machine learning techniques~\cite{Shen:2021a,Bone:2020a,Tian:2022a,Amor:2021a,Krebs:2019a,Sun:2022a,Yang:2017a,Wu:2023a,Bharati:2022a,Wu:2022a}. As we will show, our dedicated hardware-accelerated implementation~\cite{Brunn:2020a,Brunn:2021a} allows us to solve diffeomorphic image registration problems in 3 to 4 seconds on a single GPU without considering machine learning approaches.

We consider a globalized, reduced space (Gauss--)Newton--Krylov method~\cite{Mang:2015a,Mang:2019a}. For these methods to be effective, it is crucial to design a good preconditioner for solving the reduced space KKT system~\cite{Benzi:2005a}. Related work on designing preconditioners for problems similar to ours can be found in~\cite{Benzi:2011a,Simoncini:2012a,Herzog:2019a}. Another key ingredient is fast algorithms to solve the PDEs that appear in the optimality systems. In our case, the most expensive PDE operators are (hyperbolic) transport equations. We refer to~\cite{Borzi:2002a,Polzin:2016a,Hart:2009a,Benzi:2011a,Simoncini:2012a,Beg:2005a,Chen:2011a,Mang:2017a} for different numerical methods to solve these types of PDEs in the context of PDE-constrained optimization. We use a SL method~\cite{Beg:2005a,Chen:2011a,Mang:2017b,Mang:2016b}.

What separates CLAIRE~\cite{claireweb,Mang:2019a} from most existing packages for velocity-based diffeomorphic image registration, aside from the numerics, is that it features hardware-accelerated computational kernels and that it has been deployed to dedicated HPC architectures~\cite{Mang:2016b,Gholami:2017a,Mang:2019a,Brunn:2020a,Brunn:2021a}. Examples for parallel algorithms for PDE-constrained optimization problems can be found in~\cite{Akcelik:2002a,Akcelik:2006a,Biros:1999a,Biros:2005a,Biros:2005b,Biegler:2003a,Biegler:2007a,Schenk:2009a}. Surveys for parallel implementations of image registration algorithms are~\cite{Eklund:2013a,Fluck:2011a,Shackleford:2013a,Shams:2010a}. Many of these works consider low-dimensional parameterizations based on an expansion of the deformation map $\vec{y}$ in terms of smooth basis functions. Examples of GPU implementations of these approaches are~\cite{Shackleford:2010a,Modat:2010a,Shamonin:2014a}. GPU implementations of formulations similar to ours are described in~\cite{Ha:2009a,Ha:2011a,Sommer:2011a,Rehman:2009a,ValeroLara:2014a,Bone:2018a}. Our memory-distributed implementation uses MPI for parallelism and allows us to solve problems of unprecedented scale~\cite{Brunn:2020a,Mang:2016a,Gholami:2017a,Himthani:2022a}. The linear solvers and the optimizer are built on top of {\tt PETSc}/{\tt TAO}~\cite{Munson:2017a,Mills:2021a,petsc-web-page,petsc-user-ref}. Our CPU implementation~\cite{Mang:2019a,Mang:2016a,Gholami:2017a} allows us to solve problems with $\inum{3221225472}$ unknowns in \SI{2}{\minute} on 22 compute nodes (256~MPI tasks), and in less than \SI{5}{\second} if we use 342 compute nodes (\inum{4096}~MPI tasks). Our GPU implementation~\cite{Brunn:2020a,Brunn:2021a} allows us to solve clinically relevant problems (\inum{50000000} unknowns) in less than \SI{5}{\second} on a single GPU. Our multi-GPU implementation for large scale problems is described in~\cite{Brunn:2020a} and applied to large scale (biomedical) imaging data in~\cite{Himthani:2022a}. We limit the numerical results reported in this study to our GPU implementation~\cite{Brunn:2020a,Brunn:2021a}.

\subsection{Contributions}

Our contributions are as follows:
\begin{itemize}
\item We overview our work on CLAIRE---a memory distributed algorithmic framework for diffeomorphic image registration based on variational optimization problems governed by transport equations. In particular, we recapitulate our contributions presented in~\cite{Mang:2015a,Mang:2016a,Mang:2016b,Mang:2018a,Mang:2017b,Mang:2019a,Gholami:2017a,Brunn:2021b,Brunn:2020a,Brunn:2021a}.
\item We report new results and study the performance of CLAIRE for real world medical imaging data in 3D~\cite{Mang:2019a,claireweb,Brunn:2021b}. We include results for different similarity measures---normalized cross correlation and the default squared $L^2$-distance. In addition, we present results for an improved implementation in which we store the state variable and its gradient, reducing the runtime from roughly five seconds reported in prior work to slightly more than three seconds for clinically relevant problems (\inum{50000000} unknowns).
\end{itemize}

\subsection{Limitations}

CLAIRE has several limitations. First, CLAIRE only supports stationary velocity fields. Stationary velocities yield similar residuals in practical applications~\cite{Mang:2015a}. However, they are less expressive~\cite{Mang:2017a}; they only allow us to model a subset of the diffeomorphisms that can be modeled by using non-stationary $\vec{v}$. Second, CLAIRE only supports the registration of images acquired from the same modality. Implementing distance measures that allow for more complicated intensity relationships between images to be registered requires more work. Third, while our schemes for preconditioning the reduced space Hessian are effective and, in general, mesh-independent, they are not independent of the choice of the regularization parameter. Fourth, the GPU implementation only supports single precision. We trade numerical accuracy for computational throughput. For example, applying (the inverse of) high-order differential operators (e.g., biharmonic operators) results in significant numerical round-off errors and, consequently, is currently not supported. This, in conjunction with other algorithmic choices does not allow us to solve the optimization problem to arbitrary accuracy, in particular for practical parameter choices. Nonetheless, we can see that in practice we obtain an excellent agreement between the registered datasets even if the gradient of our problem is not driven to zero.

\subsection{Outline}

We present the formulation and numerical methods in~\secref{s:methods}. This includes a discussion of the mathematical framework that motivates our approach (see~\secref{s:background}), a brief recapitulation of the general problem formulation (see~\secref{s:formulation}), the optimality conditions (see~\secref{s:optcond}), and the Newton step (see~\secref{s:newtonstep}), followed by a presentation of our numerical approach (see~\secref{s:numerics}). We present some numerical results in~\secref{s:results} and conclude with~\secref{s:conclusions}.

\section{Methods}
\label{s:methods}

In the following, we present the problem formulation as well as our numerical approach and implementation aspects. We start with discussing some background material related to the considered problem formulation.

\subsection{Mathematical Foundations}\label{s:background}

Our problem formulation is related to LDDMM~\cite{Trouve:1998a,Trouve:1995a,Glaunes:2004a,Glaunes:2008a,Younes:2019a,Beg:2005a}---a mathematical framework for diffeomorphic image registration and shape matching. It builds upon the seminal work~\cite{Arnold:1966a,Arnold:1976a,Ebin:1970a}.

Let $\vec{k} = (k_1,\ldots,k_d) \in \mathbb{N}^d$, $d \in \{1,2,3\}$, denote a multi-index, and let
\[
\partial^{\vec{k}} \defeq \frac{\partial^{|\vec{k}|}}{\partial_1^{k_1}\cdots \partial_d^{k_d}}
\]

\noindent denote the differential operator of order $|\vec{k}| = \sum_{i=1}^d k_i$. Here, $\partial_i$ denotes the partial derivative with respect to the coordinate direction $x_i$ with $\vec{x} \defeq (x_1, \ldots, x_d) \in \Omega$ defined on some domain $\Omega \subseteq \mathbb{R}^d$. Moreover, let $q \in \mathbb{N}$, $1 \leq p \leq \infty$. We denote by
 \[
C^q(\Omega)
\defeq
\left\{u : \Omega \to \mathbb{R} : \partial^{\vec{k}} u \text{ is continuous for }|\vec{k}| \leq q \right\}
\]

\noindent the space of $q$-times continuously differentiable functions on $\Omega$. Moreover, let
\[
W^{q,p}(\Omega)
\defeq
\left\{
u \in L^p(\Omega) : \partial^{\vec{k}} u \in L^p(\Omega) \text{ for } 0 \leq |\vec{k}| \leq q
\right\}
\]

\noindent denote the Sobolev space with norm
\begin{equation}\label{e:sobolevnorm}
\|u\|_{q,p}
\defeq
\begin{cases}
\left(\sum_{0 \leq |\vec{k}| \leq q} \|\partial^{\vec{k}} u\|^p_p\right)^{1/p}  &\text{if } 1 \leq p < \infty,\\
\max_{0 \leq |\vec{k}| \leq q} \|\partial^{\vec{k}} u \|_{\infty} &\text{if } p = \infty.
\end{cases}
\end{equation}

\noindent Here, $\|\,\cdot\,\|_{\infty}$ denotes the standard supremum norm. Using these definitions, we denote by $C^q_0(\Omega)^d \subset C^q(\Omega)^d$ with $q \in \mathbb{N}$ the completion of the space of vector fields of class $C^q$ which along with their derivatives of order less than or equal to $q$ converge to zero at infinity. The space $C^q_0(\Omega)^d$ is a Banach space for the norm $\|u \|_{q,\infty}$. Similarly, we define the Sobolev space $W^{q,p}_0(\Omega)^d$ as a space that consists of elements with compact support on $\Omega \subseteq \mathbb{R}^d$.

We introduce a pseudo-time variable $t \in [0,1]$, a suitable Hilbert space $\fs{H}$ of smooth vector fields in $\mathbb{R}^d$, and parameterize diffeomorphisms using smooth vector fields $\vec{v} \in \fs{V}$, $\fs{V} \defeq L^r([0,1],\fs{H})$, $1 \leq r \leq \infty$, $t \mapsto \vec{v}_t \defeq \vec{v}(t,\cdot\,)$, $\vec{v}_t \in \fs{H}$. This allows us to model the flow of $\mathbb{R}^d$-diffeomorphisms $\vec{\phi}_t \defeq \vec{\phi}(t, \cdot\,)$ as the solution of the ODE
\begin{equation}\label{e:flow}
\begin{aligned}
\partial_t \vec{\phi}_t &= \vec{v}_t \circ \vec{\phi}_t && \text{for}\; t\in (0,1],\\
\vec{\phi}_t & = \operatorname{id}_{\mathbb{R}^d} && \text{for}\; t = 0,
\end{aligned}
\end{equation}

\noindent where $\operatorname{id}_{\mathbb{R}^d} : \mathbb{R}^d \to \mathbb{R}^d$, $\operatorname{id}_{\mathbb{R}^d}(\vec{x}) = \vec{x}$, is the identity transformation in $\mathbb{R}^d$ and the vector field $\vec{v}_t$ tends to zero as $\vec{x} \to \infty$; that is, we assume $\vec{v}_t \in C^q_0(\Omega)^d \supseteq \fs{H}$ for any $t \in [0,1]$. This assumption, along with suitable regularity requirements in time, guarantees that~\eqref{e:flow} admits a unique solution. Moreover, it is ensured that solutions of~\eqref{e:flow} are $\mathbb{R}^d$-diffeomorphisms~\cite{Younes:2019a}.

We assume $L^1$-integrability in time, i.e., $r = 1$~\cite{Younes:2019a}. The differentiability class $q$ and the integrability order $p$ of the Sobolev norm~\eqref{e:sobolevnorm} are chosen to stipulate adequate regularity requirements in space. A common choice for $p$ is $p = 2$. The choice of $q$ depends on the dimension $d$ of the ambient space $\Omega$. In general, we have $d\in\{1,2,3\}$. Based on the Sobolev embedding theorem~\cite{Ziemer:1989a} we observe that for $p=2$ and $q > (d/2) + 1$ the embedding $W^{q,2}_0(\Omega) \hookrightarrow C^1(\bar{\Omega})$ is compact. Since this embedding holds for all components of $\vec{v}_t$, we have that $\vec{v}_t \in \fs{H} = W^{q,2}(\Omega)^d$ with $q > 5/2$ for $d=3$ is an admissible space that yields a diffeomorphic flow $\vec{\phi}_t$ of smoothness class $1 \leq s < q - (3/2)$, $s \in \mathbb{N}$. We refer to~\cite{Dupuis:1998a,Glaunes:2008a,Younes:2019a} for a more rigorous discussion.

The set of all endpoints $\vec{y} \defeq \vec{\phi}_1$ at time $t=1$ of admissible flows $\vec{\phi}_t$ is a subgroup
\[
\fs{G} \defeq \left\{ \vec{\phi}_1 : \int_0^1 \!\!\|\vec{v}_t\|_{\fs{H}} \,\d t < \infty\right\}
\]

\noindent of $C^s$-diffeomorphisms in $\mathbb{R}^d$. This subgroup can be equipped with a right-invariant metric defined as the minimal path length of all geodesics joining two elements in $\fs{G}$~\cite{Trouve:1995a,Miller:2002a,Younes:2020a,Bauer:2014a}. The geodesic distance between $\operatorname{id}_{\mathbb{R}^d}$ and a mapping $\vec{y} \in \fs{G}$ corresponds to the square root of the kinetic energy
\begin{equation}\label{e:kinv}
\operatorname{kin}(\vec{v})
\defeq
\| \vec{v} \|_{L^2([0,1],\fs{H})}^2
=
\int_0^1 \| \vec{v}_t \|_{\fs{H}}^2 \,\d t
\end{equation}

\noindent subject to the constraint that $\vec{y}$ is equal to the solution $\vec{\phi}$ of~\eqref{e:flow} at time $t=1$ for the energy minimizing velocity $\vec{v}$. We denote this geodesic distance by $\operatorname{dist}_{\fs{G}}(\operatorname{id}_{\mathbb{R}^d},\vec{y})$,
\[
\operatorname{dist}_{\fs{G}}(\operatorname{id}_{\mathbb{R}^d},\vec{y})^2 \defeq \inf_{\vec{v} \in \fs{V}}
\left\{
\int_0^1 \!\!\| \vec{v}_t \|^2_{\fs{H}} \,\d t
: \vec{v} \in \fs{V}, \,
\vec{y} = \vec{\phi}_1,\,
\partial_t \vec{\phi}_t = \vec{v}_t\circ\vec{\phi}_t,\,
\vec{\phi}_0 = \operatorname{id}_{\mathbb{R}^d}
\right\}.
\]

\noindent The geodesic distance between two maps $\vec{y}$ and $\vec{\psi}$ is given by $\operatorname{dist}_{\fs{G}}(\operatorname{id}_{\mathbb{R}^d},\vec{y}\circ\vec{\psi}^{-1})$.

Similarly, we can measure the geodesic distance between two images $m_0,m_1 \in \fs{I}$ in terms of the kinetic energy associated with the energy minimizing $\vec{v}$ that gives rise to the diffeomorphic flow $\vec{\phi}$ that maps $m_0$ to $m_1$. To do so, we assume that the image $m_1$ is in the orbit $\fs{I}$ of the template image $m_0$ for the group $\fs{G}$ of diffeomorphism, where
\[
\fs{I} \defeq \{ m : \Omega \to \mathbb{R} : m = m_0 \circ \vec{y}^{-1}, \; \vec{y} \in \fs{G} \}.
\]

\noindent Using the geodesic distance $\operatorname{dist}_{\fs{G}}$ introduced above we have
\[
\operatorname{dist}_{\fs{I}}(m_0,m_1) \defeq \inf_{\vec{y} \in \fs{G}}
\left\{\,
\operatorname{dist}_{\fs{G}}(\operatorname{id}_{\mathbb{R}^d},\vec{y})
: m_1 = m_0 \circ \vec{y}^{-1}\right\},
\]

\noindent where $\vec{y}$ corresponds to the endpoint of the flow $\vec{\phi}$. This notion of measuring distances between deformable objects has led to the emergence of a field of study in medical image analysis referred to as \emph{computational anatomy}~\cite{Grenander:1998a,Miller:2004a,Miller:2002a,Younes:2009a,Miller:2015a}.

Putting everything together, we can formulate the diffeomorphic matching of the template image $m_0$ to the reference image $m_1$ as a variational optimization problem. We stated initially that we seek a diffeomorphic map $\vec{y}$ such that $m_0 \circ \vec{y}^{-1} = m_1$. This is an ill-posed problem; we try to estimate a vector field given scalar data. Consequently, a solution $\vec{y}$ may not exist, and if it exists, it may not be unique or depend continuously on the data. To alleviate the ill-posedness we introduce a regularization model that rules out unwanted solutions. For example, we can restrict ourselves to maps $\vec{y}$ that are close to the identity $\operatorname{id}_{\mathbb{R}^d}$, i.e., we penalize the distance between $\operatorname{id}_{\mathbb{R}^d}$ and $\vec{y}$. To alleviate existence issues, we relax the exact matching requirement to $m_0 \circ \vec{y}^{-1} \approx m_1$. We do so by introducing a distance that measures the proximity between the deformed template image $m_0 \circ \vec{y}^{-1}$ and the reference image $m_1$. In conclusion, we seek $\vec{y} \in \fs{G} \subset \operatorname{diff}(\mathbb{R}^d)$ as a solution to
\[
\minopt_{\vec{y} \in \fs{G}} \;\;
\operatorname{dist}(m_0 \circ \vec{y}^{-1},m_1)
+ \frac{\alpha}{2}\operatorname{dist}_{\fs{G}}(\vec{y}, \operatorname{id}_{\mathbb{R}^d}).
\]

\noindent We can reformulate the variational problem above as an optimal control problem governed by~\eqref{e:flow}~\cite{Beg:2005a,Dupuis:1998a}. We have
\begin{subequations}\label{e:varoptlddmm}
\begin{align}
\minopt_{\vec{\phi} \in \fs{F}_{\text{ad}},\;\vec{v}\in\fs{V}_{\text{ad}}} \;\;
& \operatorname{dist}( m_0 \circ \vec{\phi}_1^{-1},m_1)
+ \frac{\alpha}{2}\operatorname{kin}(\vec{v})\\
\begin{aligned}
\text{subject to} \;\;\\{}\\
\end{aligned}
&
\begin{aligned}
\partial_t \vec{\phi}_t & = \vec{v}_t \circ \vec{\phi}_t && \text{in } (0,1],\\
\vec{\phi}_t & = \operatorname{id}_{\mathbb{R}^d} && \text{for}\; t = 0,
\end{aligned}
\end{align}
\end{subequations}

\noindent where the first term in the objective functional measures the discrepancy between the deformed template image $m_0 \circ \vec{\phi}_1^{-1}$ and the reference image $m_1$, the second term denotes the kinetic energy in~\eqref{e:kinv} and the parameter $\alpha > 0$ balances their contribution. The norm $\|\vec{v}_t \|^2_{\fs{H}}$ in the definition of the kinetic energy in~\eqref{e:kinv} is typically modelled as
\[
\|\vec{v}\|_{\fs{H}}^2
= \langle\vec{v},\vec{v}\rangle_{\fs{H}}
= \langle\op{B}\vec{v},\op{B}\vec{v}\rangle_{\mathbb{R}^d},
= \langle\op{L}\vec{v},\vec{v}\rangle_{\mathbb{R}^d},
\]

\noindent where $\op{L} : \fs{V} \to \fs{V}^\dual$, $\op{L} = \op{B}^\ast\op{B}$, is a differential operator of adequate order. A common choice for $\op{B}$ is a symmetric, positive definite Helmholtz operator of the form $\op{B} \defeq (\beta \operatorname{id} - \ilap_d)^\gamma$, $\beta, \gamma > 0$~\cite{Beg:2005a}, where $\ilap_d \vec{u}(\vec{x}) \defeq ( \ilap u_1(\vec{x}), \ldots, \ilap u_d(\vec{x}))$, $\ilap \defeq \sum_{i=1}^d \partial_{ii}$ for any $\vec{u} : \bar{\Omega} \to \mathbb{R}^d$.

Other data structures than images that can be registered within this framework are landmarks~\cite{Glaunes:2004b,Joshi:2000a}, curves~\cite{Durrleman:2010a,Durrleman:2008a}, surfaces~\cite{Zhang:2021a,Azencott:2010a,Hartman:2023a,Mang:2023a,Durrleman:2010a,Glaunes:2006a,Kurtek:2012a,Arguillere:2016a}, tensor fields~\cite{Cao:2006a} or functional data on manifolds. We refer to~\cite{Beg:2005a,Azencott:2010a,Zhang:2021a,Mang:2023a,Polzin:2020a} for numerical methods to solve the control problem in \eqref{e:varoptlddmm}.

\subsection{Variational Problem Formulation}\label{s:formulation}

In this section, we review the problem formulation considered in CLAIRE. We assume that $m_i : \Omega \to \mathbb{R}$, $i=0,1$, are smooth $C^1$-functions compactly supported on $\Omega$. As stated in \secref{s:intro}, we formulate diffeomorphic image registration as a PDE-constrained optimization problem of the general form~\eqref{e:varopt}. This is different from the ODE-constrained optimization problem~\eqref{e:varoptlddmm}. Motivated by the formulation discussed above, we introduce a pseudo-time variable $t \in [0,1]$ and invert for a smooth, time-dependent velocity field $\vec{v} \in \fs{V}$~\cite{Mang:2015a}. However, to reduce the computational complexity, our hardware-accelerated implementation no longer inverts for a time-dependent velocity $\vec{v}$ but for a \emph{stationary velocity field} $\vec{v} : \Omega \to \mathbb{R}^d$. This not only reduces the complexity of the optimization problem but also simplifies the implementation. We note that stationary velocities no longer define a Riemannian metric as described in \secref{s:background}. However, we still generate diffeomorphic transformations. Moreover, we did not observe a deterioration in registration accuracy when comparing results to a non-stationary implementation~\cite{Mang:2015a}. Related work by other groups that use stationary $\vec{v}$ can be found in~\cite{Arsigny:2006a,Ashburner:2007a,Hernandez:2009a,Lorenzi:2013a,Lorenzi:2013b,Vercauteren:2009a}.

In its simplest form, the PDE constraint $c$ in~\eqref{e:varopt} for stationary $\vec{v}$ is given by the \emph{hyperbolic transport equation}
\begin{subequations}
\label{e:transport}
\begin{align}
\partial_t m(t,\vec{x}) + \igrad m(t,\vec{x}) \cdot \vec{v}(\vec{x})  & = 0
&& \text{in}\,\,(0,1] \times \Omega,
\\
m(t,\vec{x}) &=  m_0(\vec{x})
&& \text{in}\,\,\{0\} \times \Omega.
\end{align}
\end{subequations}

\noindent The solution of this PDE is the transported intensities $m$ of the template image $m_0$. The endpoint $m(1) \defeq m(1, \,\cdot\,)$ at time $t=1$ corresponds to the deformed template image for some trial velocity $\vec{v}$.

The second building block of our variational problem formulation is the distance functional $\dist : \fs{I} \times \fs{I} \to \mathbb{R}$ in~\eqref{e:varopt:objective} that quantifies the discrepancy between the transported template image $m(1) \defeq m(\,\cdot\,, 1)$ at time $t=1$ and the reference image $m_1$. A common choice for this terminal (endpoint) cost is given by the squared $L^2$-distance
\[
\dist(m(1), m_1)
= \frac{1}{2}\int_{\Omega} ( m(1,\vec{x}) - m_1(\vec{x}) )^2 \, \d \vec{x}.
\]

\noindent While this is a common choice in many diffeomorphic image registration packages~\cite{Beg:2005a,Mang:2019a}, this distance measure can only be used for registering images acquired using the same imaging modality. We present an alternative in the appendix.

The last building block is the regularization functional $\reg : \fs{V} \to \mathbb{R}$. Motivated by the problem formulation presented in \secref{s:background}, we use
\[
\reg(\vec{v})
= \frac{\alpha}{2} \|\vec{v}\|_{\fs{H}}^2
= \frac{\alpha}{2} \langle\op{L}\vec{v},\vec{v}\rangle_{\mathbb{R}^d},
\]

\noindent where $\op{L} : \fs{V} \to \fs{V}^\dual$ is a differential operator of adequate order. CLAIRE, in general, features $H^1$-, $H^2$- and $H^3$-norms and semi-norms for the regularization of $\vec{v}$~\cite{Mang:2015a,Mang:2019a,Mang:2016a}. The default regularization operator is an $H^1$-seminorm, i.e., $\op{L} = -\ilap_d$, with an additional $H^1$-norm that penalizes the divergence of the velocity~\cite{Mang:2019a,Brunn:2020a,Brunn:2021a,Himthani:2022a}. We provide additional details in the appendix.

Putting everything together, we arrive at the PDE-constrained optimization problem
\begin{subequations}\label{e:varopt-claire}
\begin{align}
\minopt_{m \in \fs{M}_{\text{ad}}, \, \vec{v} \in \fs{V}_{\text{ad}}}\quad
&
\frac{1}{2}\int_{\Omega} (m(1,\vec{x}) - m_1(\vec{x}) )^2\, \d \vec{x}
+ \frac{\alpha}{2} \|\vec{v}\|_{\fs{H}}^2
\\
\begin{aligned}
\st \\ \\
\end{aligned}
\quad&
\begin{aligned}
\partial_t m(t,\vec{x}) + \igrad m(t,\vec{x}) \cdot \vec{v}(\vec{x})  & = 0
&& \text{in}\,\,(0,1] \times \Omega,
\\
m(t,\vec{x}) &=  m_0(\vec{x})
&& \text{in}\,\,\{0\} \times \Omega.
\end{aligned}
\end{align}
\end{subequations}

\noindent Similar problem formulations have been considered in~\cite{Borzi:2002a,Borzi:2002b,Chen:2012a,Hart:2009a}. For simplicity, we discuss the numerical methods for the problem formulation in~\eqref{e:varopt-claire}. However, we note that we considered different variants of this formulation in our past work~\cite{Mang:2015a,Mang:2016a,Mang:2019a,Brunn:2020a}. We discuss these in greater detail in the appendix.

\subsection{Optimality Conditions}\label{s:optcond}

In the present work, we consider an \emph{optimize-then-discretize} approach. The advantages of this approach are that the formal optimality conditions are straightforward to derive. They also retain interpretability; for example, we will see that the adjoint equation of the transport equation in~\eqref{e:transport} represents a continuity equation for the image mismatch (see~\eqref{e:adjoint}). Moreover, one can freely decide on the numerical methods to solve the PDEs associated with the optimality conditions. A disadvantage of this approach is that the discrete gradient is (in general) not consistent with the discretized objective functional (contingent on the numerical scheme used for discretization). Consequently, it is not possible to solve the variational optimization problem with arbitrary accuracy (i.e., to machine precision). A \emph{discretize-then-optimize} approach guarantees that the discrete gradient is consistent with the discretized objective functional. However, depending on the discretization this approach also has drawbacks. We refer to~\cite{Gunzburger:2003a} for a general discussion and to~\cite{Mang:2017a,Polzin:2020a} for examples of \emph{discretize-then-optimize} implementations for a problem of the form~\eqref{e:varopt-claire}.

We consider the method of Lagrange multipliers to solve~\eqref{e:varopt}. We introduce the dual variable $\lambda : [0,1] \times \bar{\Omega} \to \mathbb{R}$, $\lambda \in \fs{M}^\dual$, for the transport equation \eqref{e:transport}. The Lagrangian functional is given by
\begin{equation}\label{e:lagrangian}
\begin{aligned}
\ell(\vec{\Xi}) & =
\frac{1}{2}\int_{\Omega} ( m(1,\vec{x}) - m_1(\vec{x}) )^2 \, \d \vec{x}
+ \frac{\alpha}{2} \langle\op{L}\vec{v},\vec{v}\rangle_{\mathbb{R}^d}
\\
& \quad
+ \int_0^1\langle\partial_t m  + \igrad m \cdot \vec{v}, \lambda \rangle_{L^2(\Omega)} \d t
+ \langle m(0) - m_0, \lambda \rangle_{L^2(\Omega)}
\end{aligned}
\end{equation}

\noindent where $\vec{\Xi} \defeq (\vec{v},m,\lambda) \in \fs{V} \times \fs{M} \times \fs{M}^\dual$.

Computing first variations with respect to the control variable $\vec{v}$ yields the reduced gradient
\begin{equation}\label{e:rgrad}
\vec{g}[\vec{v}](\vec{x})
\defeq \alpha \op{L} [\vec{v}](\vec{x})
+ \int_0^1 \lambda(t,\vec{x}) \igrad m(t,\vec{x}) \, \d t.
\end{equation}

To be able to evaluate the reduced gradient we require the state variable $m \in \fs{M}$ and the dual variable $\lambda \in \fs{M}^\dual$. We can find the state variable by solving~\eqref{e:transport} forward in time. Formally, this equation is obtained by computing the first variations of $\ell$ in~\eqref{e:lagrangian} with respect to $\lambda$. The dual variable $\lambda$ can be found by solving the adjoint equations backward in time. Formally, the adjoint equations are found by computing variations of $\ell$ in~\eqref{e:lagrangian} with respect to $m$. We obtain
\begin{subequations}
\label{e:adjoint}
\begin{align}
-\partial_t \lambda(t,\vec{x}) - \idiv \lambda(t,\vec{x}) \vec{v}(\vec{x})  & = 0
&& \text{in}\,\,[0,1) \times \Omega,
\\
\lambda(t,\vec{x}) &=  -(m(1,\vec{x}) - m_1(\vec{x}))
&& \text{in}\,\,\{1\} \times \Omega,
\label{e:adjfinal}
\end{align}
\end{subequations}

\noindent subject to periodic boundary conditions on $\partial\Omega$. Notice that this equation represents a continuity equation; we transport the mismatch between the deformed template image $m(1,\,\cdot\,)$ and the reference image $m_1$ backward in time. If we change the distance measure $\operatorname{dist} : \fs{I} \times \fs{I} \to \mathbb{R}$ in \eqref{e:varopt-claire}, the final conditions in~\eqref{e:adjfinal} will change.

\subsection{Newton Step}\label{s:newtonstep}

We consider a (Gauss--)Newton--Krylov method for numerical optimization~\cite{Mang:2015a}. We provide more details in~\secref{s:numerics}. The PDE operators associated with the Hessian can be found by formally computing second-order variations of the Lagrangian $\ell$ in~\eqref{e:lagrangian}. The expression for the Hessian matvec---i.e., the application of the Hessian to a vector $\vec{\tilde{v}}$---is given by
\begin{equation}\label{e:matvec}
\begin{aligned}
\op{H}[\vec{\tilde{v}}](\vec{x})
& =
\op{H}_{\text{reg}}[\vec{\tilde{v}}](\vec{x})
+ \op{H}_{\text{dat}}[\vec{\tilde{v}}](\vec{x})
\\
& =
\alpha\op{L}[\vec{\tilde{v}}](\vec{x})
+ \int_0^1
\left\{ \tilde{\lambda}(t,\vec{x}) \igrad m(t,\vec{x})
+ \lambda(t,\vec{x}) \igrad \tilde{m}(t,\vec{x})\right\}
\, \d t.
\end{aligned}
\end{equation}

\noindent The variable $\vec{\tilde{v}} : \bar{\Omega} \to \mathbb{R}^d$, $\vec{v} \in \fs{V}$, represents the incremental control variable, i.e., the search direction for $\vec{v}$. The operators $\op{H}_{\text{reg}}$ and $\op{H}_{\text{dat}}$ denote the regularization part and the data part of the reduced space Hessian, respectively. For the latter, the dependence on $\vec{\tilde{v}}$ is hidden in the incremental PDE operators. Given a candidate $\vec{v}$ and a candidate $\vec{\tilde{v}}$ we require the state variable $m : [0,1] \times \bar{\Omega} \to \mathbb{R}$, the dual variable $\lambda : [0,1] \times \bar{\Omega} \to \mathbb{R}$, the incremental state variable $\tilde{m} : [0,1] \times \bar{\Omega} \to \mathbb{R}$, and the incremental adjoint variable $\tilde{\lambda} : [0,1] \times \bar{\Omega} \to \mathbb{R}$. We can find the state and dual variables during the evaluation of the reduced gradient in~\eqref{e:rgrad}. The incremental state variable can be found by solving
\begin{subequations}\label{e:incstate}
\begin{align}
\partial_t \tilde{m}(t,\vec{x})
+ \igrad \tilde{m}(t,\vec{x}) \cdot \vec{v}(t,\vec{x})
+ \igrad m(t,\vec{x}) \cdot \vec{\tilde{v}}(t,\vec{x})
& = 0
&& \text{in}\,\,(0,1] \times \Omega,
\\
\tilde{m}(t,\vec{x}) &= 0
&& \text{in}\,\,\{0\} \times \Omega,
\end{align}
\end{subequations}

\noindent subject to periodic boundary conditions on $\partial\Omega$ forward in time. We can find the incremental dual variable $\tilde{\lambda}$ by solving
\begin{subequations}\label{e:incadj}
\begin{align}
-\partial_t \tilde{\lambda}(t,\vec{x})
+ \idiv(\tilde{\lambda}(t,\vec{x}) \vec{v}(\vec{x}) + \lambda(t,\vec{x}) \vec{\tilde{v}}(\vec{x})) & = 0
&& \text{in}\,\,[0,1) \times \Omega,
\\
\tilde{\lambda}(t,\vec{x}) &= -\tilde{m}(t,\vec{x})
&& \text{in}\,\,\{1\} \times \Omega,
\end{align}
\end{subequations}

\noindent subject to periodic boundary conditions on $\partial\Omega$ backward in time. Consequently, each time we apply $\op{H}$ to a vector, we have to solve two PDEs.

\subsection{Numerics}\label{s:numerics}

The numerical implementation discussed below is based on the computational kernels described in~\cite{Mang:2015a,Mang:2019a,Mang:2017b,Mang:2016b,Gholami:2017a,Brunn:2020a,Brunn:2021a,Mang:2019a}. The hardware-accelerated CPU implementation is described in~\cite{Mang:2016b,Mang:2019a,Gholami:2017a}. The GPU implementation is described in~\cite{Brunn:2020a,Brunn:2021a}.

We note that the GPU implementation is only available in single precision. This poses several numerical challenges. In particular, we observed that our scheme does not allow us to solve the optimization problem to arbitrary accuracy. This is due to the accumulation of numerical errors, dominated by the time integration and numerical differentiation. Moreover, the numerical gradient is inconsistent with the objective function. This is caused by our particular choice of a numerical time integrator in conjunction with an optimize-then-discretize approach.

\subsubsection{Discretization}

We consider a nodal discretization in space. That is, we subdivide the spatial interval $\Omega = [-\pi,\pi)^d \subset \mathbb{R}^d$ into $n_i \in \mathbb{N}$ cells of width $h_i > 0$, $i=1,\ldots,d$, along each spatial direction $x_i$, $i=1,\ldots,d$. The width of the cells $\vec{h} = (h_1, \ldots, h_d) \in \mathbb{R}^d$ along each spatial direction is given by $h_i = 2 \pi / n_i$, $i = 1, \ldots, d$. Let $\vec{x}_{\vec{l}} \in \mathbb{R}^d$ denote a mesh point at index $\vec{l} = (l_1,\ldots,l_d) \in \mathbb{N}^d$, $1 \leq l_i \leq n_i$. The coordinates $\vec{x}_{\vec{l}}$ are computed according to \[\vec{x}_{\vec{l}} = ((\vec{n}/2) - \vec{l}) \odot \vec{h}.\] Here, $\odot$ denotes an elementwise multiplication between two vectors (Hadamard product) and $\vec{n} = (n_1, \ldots, n_d) \in \mathbb{N}^d$ represents the number of mesh points along each spatial direction. We denote the resulting mesh by $\di{\Omega} = (\di{x}_{\vec{l}}) \in \mathbb{R}^{d, n_1, \ldots, n_d}$. Similarly, we subdivide the unit time interval $[0,1]$ into a uniform mesh with step size $h_t = 1 / n_t$. We discretize integrals using a trapezoidal rule.

\subsubsection{Time Integration} We use a SL method for numerical time integration~\cite{Staniforth:1991a}. The prototype implementation of this time integrator is described in~\cite{Mang:2017b}. Different variants of hardware-accelerated implementations are described in~\cite{Mang:2016b,Brunn:2020a,Brunn:2021a}. Other works that consider a SL scheme in a similar context are~\cite{Franccois:2021a,Beg:2005a}.

The SL scheme is a hybrid between Eulerian and Lagrangian methods. It is unconditionally stable, i.e., we can select the time step size $h_t > 0$ solely based on accuracy considerations. To apply the SL scheme to the PDEs that appear in our optimality system we need to bring them into the general form
\begin{equation}\label{e:matderiv}
\d_t u = f(t,u, \vec{v}, \ldots)
\end{equation}

\noindent where $u : [0,1] \times \Omega \to \mathbb{R}$ denotes an arbitrary scalar function, $\d_t \defeq \partial_t + \vec{v} \cdot \igrad$ denotes the \emph{material derivative}, and the right-hand-side $f$ represents all remaining terms. To obtain this representation for the equations considered here, we use the vector calculus identity $\idiv u\vec{v} = u \idiv \vec{v} + \igrad u \cdot \vec{v}$.

\begin{figure}
\centering
\includegraphics[height=4cm]{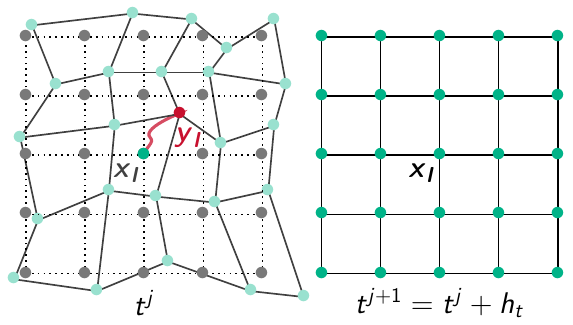}
\caption{Illustration of the computation of the characteristic $\vec{y}$ in the SL scheme. In the SL scheme, we have to compute the departure points at time $t^j$. To do so, we start with a regular grid at time $t^{j+1}$ and solve for the characteristic $\vec{y}_{\vec{l}}$ at a given point $\vec{x}_{\vec{l}}$ at mesh index $\vec{l}$ backward in time (green line in the graphic on the left). The deformed grid configuration is overlaid onto the initial regular grid at time $t^j$. (Figure modified from \cite{Brunn:2020a,Mang:2017a}.)}
\label{f:illustration-sml}
\end{figure}

In the first step, we have to compute the characteristic $\vec{y} : [t^j,t^{j+1}] \to \mathbb{R}^d$ along which particles flow between the timepoint $t^j$ and $t^{j+1}$, $j = 1, \ldots, n_t$. The question we seek to answer is where particles at time $t^{j+1}$ originate from given data at time $t^j$. That is, we compute the departure point $\vec{y}$ at time $t = t^j$. To compute this characteristic, we solve the ODE $\d_t \vec{y} = \vec{v} \circ \vec{y}$ for $t \in [t^j, t^{j+1})$ with $\vec{y} = \vec{x}$ for $t = t^{j+1}$ backward in time.  We illustrate the computation of the departure point in \figref{f:illustration-sml}.

In our implementation, we compute the characteristics $\vec{y}$ using an RK2 method. Notice that the velocity $\vec{v}$ is constant in time; this simplifies the computation of the transported quantities considerably. Let $\di{v}_{\vec{l}} \defeq \di{v}(\di{x}_{\vec{l}})$ denote the discretized velocity at a given mesh point $\di{x}_{\vec{l}} \in \di{\Omega}$. We obtain the $\vec{l}$-th query point $\di{y}_{\vec{l}} \in \mathbb{R}^d$ associated with $\di{x}_{\vec{l}} \in \di{\Omega}$ according to
\[
\begin{aligned}
\di{\tilde{y}}_{\vec{l}} & \gets \di{x}_{\vec{l}} - h_t\di{v}(\di{x}_{\vec{l}})\\
\di{y}_{\vec{l}} & \gets \di{x}_{\vec{l}} - \frac{h_t}{2}\left(\di{v}(\di{x}_{\vec{l}}) + \di{v}(\di{\tilde{y}}_{\vec{l}})\right).
\end{aligned}
\]

\noindent The intermediate query points $\di{\tilde{y}}_{\vec{l}}$ and the final query point $\di{y}_{\vec{l}}$ (i.e., the departure point), are---in general---off-grid locations. Therefore, evaluating quantities of interest at these locations requires interpolation (see \secref{s:interpolation} for details). If~\eqref{e:matderiv} is homogeneous, i.e., $f = 0$, we only interpolate the transported quantity to obtain its value at the departure point $\di{y}_{\vec{l}}$ at time $t^j$ and assign the resulting value to the regular mesh point $\di{x}_{\vec{l}} \in \di{\Omega}$ at time $t^{j+1}$. That is,
\[
u(t^{j+1},\di{x}_{\vec{l}}) \gets u(t^{j},\di{y}_{\vec{l}}).
\]

If~\eqref{e:matderiv} is not homogeneous, i.e., $f \not = 0$, we have to solve the ODE~\eqref{e:matderiv} along the characteristic $\vec{y}$ forward in time. We do so using an RK2 scheme. That is,
\[
\begin{aligned}
f_0 & \gets f(t^j,u(t^j,\di{y}_{\vec{l}}), \di{v}(\di{y}_{\vec{l}}), \ldots) \\
\tilde{u}(t^{j+1},\di{x}_{\vec{l}}) & \gets u(t^j,\di{y}_{\vec{l}}) + h_t f_0 \\
f_1 & \gets f(t^{j+1},\tilde{u}(t^{j+1},\di{x}_{\vec{l}}), \di{v}(\di{x}_{\vec{l}}), \ldots) \\
u(t^{j+1},\di{x}_{\vec{l}}) & \gets u(t^j,\di{y}_{\vec{l}}) - h_t(f_0 + f_1) / 2.
\end{aligned}
\]

\noindent Again, quantities evaluated at the query point $\di{y}_{\vec{l}}$ at time $t^j$ need to be interpolated. We note that these functions live on a curvilinear mesh (see \figref{f:illustration-sml}). Since we use spectral methods with a Fourier basis we cannot evaluate the differential operators that appear in $f$ on such a mesh. As a remedy, we do not compute the derivative on this curvelinear mesh but interpolate (i.e., transport) the derivatives evaluated on a regular mesh instead.

\subsubsection{Differentiation}

In our past work, we considered pseudo-spectral methods with a Fourier basis for numerical differentiation~\cite{Mang:2015a,Mang:2017b,Mang:2016b,Mang:2019a}.

\begin{figure}
\centering
\includegraphics[width=\textwidth]{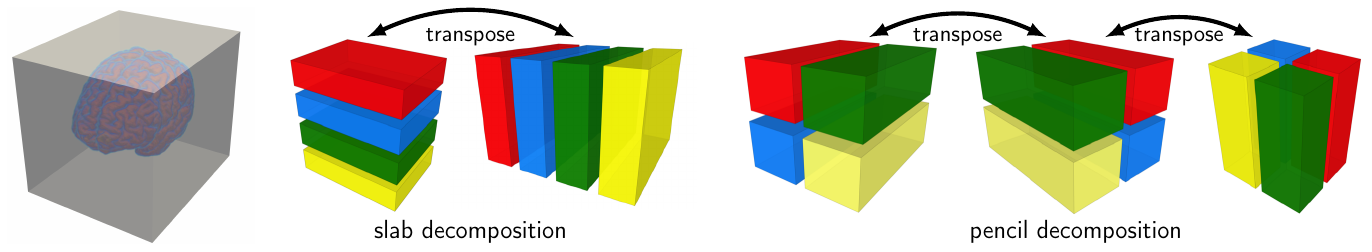}
\caption{Domain decomposition for memory-distributed implementation. In each case, we assume that we use four MPI tasks to distribute our data (e.g., four GPUs or four nodes). Left: 3D volume rendering of medical imaging data set (brain image). Middle: Slab decomposition (1D domain decomposition) considered in our GPU implementation~\cite{Brunn:2020a}. We decompose the spatial domain in the outer-most dimension. We transpose the data only once. On the right we illustrate the pencil decomposition (2D domain decomposition) of the data considered in our CPU implementation~\cite{Mang:2016b,Mang:2019a,Gholami:2017a}. The computation in this data layout involves three transposes.\label{f:decomposition}}
\end{figure}

For our GPU implementation~\cite{Brunn:2020a,Brunn:2021a} we have designed a mixed-precision approach to improve scalability and computational throughput. We consider $8^{\text{th}}$-order finite differences for first-order derivatives (i.e., the gradient and divergence operators). Higher order derivative operators (e.g., the Laplacian operator $\ilap$) and their inverse are implemented using a pseudo-spectral discretization with a Fourier basis. That is, we model an arbitrary function $u : \bar{\Omega} \to \mathbb{R}$ discretized on a regular mesh $\di{\Omega}$ at grid points $\di{x}_{\vec{l}}\in \di{\Omega}$, $\vec{l} = (l_1, \ldots, l_d) \in \mathbb{N}^d$, $l_i = 1,\ldots, n_i$, as $u_{\vec{l}} \defeq u(\di{x}_{\vec{l}})$,
\[
u_{\vec{l}} =
\sum_{k_1= -(n_1/2) + 1}^{n_1/2} \cdots\sum_{k_d= -(n_d/2) + 1}^{n_d/2}
\hat{u}_{\vec{k}} \exp(- \text{i} \langle \vec{k}, \di{x}_{\vec{l}}\rangle_{\mathbb{R}^d})
\]

\noindent with $\vec{k} \in \mathbb{Z}^d$ and spectral cofficients $\hat{u}_{\vec{k}}$. This spectral representation is the reason why we assume periodic boundary conditions in our continuous model. We note that images may not necessarily be periodic functions. We can address this by zero-padding the datasets and applying a mollifier close to the boundary $\partial \Omega$. The mapping between the coefficients $\{u_{\vec{l}}\}$ and $\{\hat{u}_{\vec{k}}\}$ is done using forward and inverse FFTs. In our CPU implementation, we considered a pencil decomposition~\cite{Mang:2016b,Mang:2019a,Gholami:2017a} (see \figref{f:decomposition}; right). Here, 1D FFTs along each spatial direction are computed based on the {\tt FFTW} library. FFTs along other directions are then obtained by transposing the data, resulting in large communication costs. For the single GPU implementation described in \cite{Brunn:2021a} we switched to \texttt{cuFFT} for 3D FFTs. The multi-GPU implementation described in~\cite{Brunn:2020a} uses a combination of \texttt{cuFFT} and a new 2D slab decomposition (see \figref{f:decomposition}; middle). This enables us to utilize the highly optimized 2D \texttt{cuFFT} on each GPU. We decompose the spatial domain in the outer-most dimension (i.e., $x_1$) and in the spectral domain in the $x_2$ direction. Consequently, the inner-most $x_3$ direction remains continuous in memory. This reduces misaligned memory access for the communication of the transpose operations.  The real-to-complex transformation is divided into three steps: First,  we execute \texttt{cuFFT}'s batched 2D FFTs in the plane spanned by the $x_2$ and $x_3$ axis. Then, we transpose the complex data to a decomposition in $x_2$ direction. Then, we apply \texttt{cuFFT}'s batched 1D FFTs to the $x_1$ direction, which is non-continuous in memory. For the inverse complex-to-real transformation, these three steps are executed in reverse order, using the respective inverse transformations. For the execution on multiple GPUs, we use CUDA-aware MPI to eliminate expensive on-node host-device transfers~\cite{Brunn:2020a}.

The $8^{\text{th}}$-order finite difference approximation of the first order derivative along the $i$th coordinate direction $x_i$ at a mesh point $\di{x}_{\vec{l}} \in \di{\Omega}$ is given by
\[
\partial_i u(\di{x}_{\vec{l}}) \approx
\frac{u^-_i(\di{x}_{\vec{l}})  +  u^+_i(\di{x}_{\vec{l}})}{840 h_i}
\]

\noindent with $u^-_i(\di{x}_{\vec{l}}) \defeq  3u_{\vec{l} - 4 \vec{e}_i} - 32u_{\vec{l} - 3 \vec{e}_i} + 168u_{\vec{l} - 2 \vec{e}_i} - 672 u_{\vec{l} - \vec{e}_i}$,  $u^+_i(\di{x}_{\vec{l}}) \defeq  672 u_{\vec{l} + \vec{e}_i} - 168u_{\vec{l} + 2 \vec{e}_i} + 32u_{\vec{l} + 3 \vec{e}_i} - 3u_{\vec{l} + 4 \vec{e}_i}$, and unit vectors $\vec{e}_i \in \{0,1\}^d$, where the $j$-th entry is one for $j = i$ and zero otherwise.

The key limiting factor to obtain optimal strong and weak scalability for our method is the communication costs associated with the FFT. We refer to \cite{Brunn:2020a,Brunn:2021a} for additional details.

\subsubsection{Interpolation}\label{s:interpolation}

In our past work~\cite{Brunn:2021a,Brunn:2020a}, we implemented different interpolation models. To obtain the interpolated value of a function $u$ at an aribrary query point $\di{x} = (x_1,\ldots,x_d)\in\mathbb{R}^d$ we evaluate
\[
u(\di{x})
=
\sum_{k_1 = 1}^{p+1}
\ldots
\sum_{k_d = 1}^{p+1}
c_{k_1,\ldots,k_d}\prod_{i=1}^d \phi_{k_i}(x_i),
\]

\noindent where $p \in \mathbb{N}$ denotes the polynomial order, $\phi_j : \mathbb{R} \to \mathbb{R}$, $j= 1, \ldots,p+1$, are the polynomial basis functions, and $c_{k_1,\ldots,k_d} \in \mathbb{R}$ are the coefficients. For Lagrange interpolation, the coefficients are identical to the grid values of the discretized function $u$ and $\phi_j$ are the Lagrange polynomials.

The numerical accuracy and computational performance of different variants of the interpolation kernel are discussed in~\cite{Brunn:2021a}. The multi-GPU implementation is described in~\cite{Brunn:2020a}. We use NVIDIA's libraries for texture-based trilinear interpolation~\cite{Sigg:2005a}. In~\cite{Brunn:2021a}, we also transferred our CPU kernels for cubic Lagrange interpolation~\cite{Mang:2016a,Gholami:2017a,Mang:2019a} to the GPU. We implemented two variants, one that uses texture units for lookup tables and one that implements texture-based interpolation~\cite{Brunn:2021a}. The latter implementation is similar to~\cite{Ruijters:2008a}; it yields higher computational throughput at lower accuracy. Lastly, we have developed a texture-based B-spline interpolation, the computational kernels of which are inspired by~\cite{Ruijters:2008a,Guijters:2012a,Champagnat:2012a}. For the execution on multiple GPUs, we use CUDA-aware MPI to eliminate expensive on-node host-device transfers~\cite{Brunn:2020a}. We utilize the \texttt{thrust} library~\cite{thrustweb} to determine which query points need to be processed by which GPU, thereby completely eliminating host-side computation. We use sparse point-to-point communication to send points to other processors. We adaptively allocate memory for the respective MPI buffers. We do this by computing an estimate of the maximal displacement of grid points along the computed trajectories based on the CFL number of the velocity field. We refer to~\cite{Brunn:2020a,Brunn:2021a} for additional details.

\subsubsection{Optimization}\label{s:optimization}

We use an iterative method globalized by an Armijo line search~\cite{Nocedal:2006a,Boyd:2004a}. The outer iterations of our algorithm are summarized in \algref{a:outer-iteration} in the appendix. At (outer) iteration $k \in \mathbb{N}$, we update the iterate $\di{v}^{(k)} $ according to
\[
\di{B}^{(k)} \di{\tilde{v}}^{(k)} = -\di{g}^{(k)}, \quad
\di{v}^{(k+1)} = \di{v}^{(k)} + \gamma^{(k)} \di{\tilde{v}}^{(k)},
\]

\noindent where $\di{\tilde{v}}^{(k)} \in \mathbb{R}^{dn}$ denotes the search direction, $\di{B}^{(k)} \in \mathbb{R}^{dn,dn}$ is a positive-definite matrix, and $\gamma^{(k)} > 0$ is the step size.

For $\di{B}^{(k)} = \operatorname{diag}(1,\ldots,1) \in \mathbb{R}^{dn,dn}$ the scheme above corresponds to a gradient descent algorithm. In~\cite{Mang:2015a}, we consider a preconditioned gradient descent algorithm. This scheme is more stable and yields an improved convergence behavior. The preconditioner is the inverse of the regularization operator $\mathcal{L}$. That is, $\di{B}^{(k)} = \alpha\di{L}$, where $\di{L} \succeq \di{0}$ denotes the discretization of $\mathcal{L}$. This scheme can be viewed as a Picard iteration. We note that our spectral discretization allows us to apply the inverse of this operator with vanishing costs; the complexity of inverting $\di{L}$ is $\mathcal{O}(n \log n)$ regardless of the Sobolev norm we consider. If the operator $\mathcal{L}$ has a non-trivial kernel, we set the spectral coefficients that are zero to one before inverting $\alpha\di{L}$. Consequently, the search direction $\di{\tilde{v}}^{(k)}$ is given by $\di{\tilde{v}}^{(k)} = -(\alpha\di{L})^{-1} \di{g}^{(k)}$.

In addition, we have designed a (Gauss)--Newton--Krylov algorithm for numerical optimization~\cite{Mang:2015a,Mang:2019a}. Here, $\di{B}^{(k)}$ corresponds to the Hessian matrix $\di{H}^{(k)}$ at (outer) iteration $k$. Consequently, we have to invert $\di{H}^{(k)}$ at each iteration. We note that forming and storing $\di{H}^{(k)}$ results in prohibitive computational costs and memory requirements. As such, we cannot use direct methods~\cite{Duff:2017a,Davis:2006a}. Instead, we use iterative methods to invert $\di{H}^{(k)}$. In particular, we use matrix-free Krylov subspace methods---more precisely, a PCG algorithm~\cite{Hestenes:1952a}---to compute the action of the inverse of $\di{H}^{(k)}$ on the vector $-\di{g}^{(k)}$. As such, we only require an expression for the Hessian matvec. This is precisely what is given by~\eqref{e:matvec}. Thus, we need to evaluate~\eqref{e:matvec} at every inner iteration of our Krylov-subspace method. This involves solving the PDEs~\eqref{e:incstate} and~\eqref{e:incadj} at every inner iteration of the PCG algorithm. These matvecs constitute the most expensive part of our algorithm. We summarize this algorithm in \algref{a:inner-iteration} in the appendix.

We note that we can use other iterative methods to compute the action of the inverse of the Hessian. In fact, we have tested different methods. In our experiments, we did not observe any issues with the PCG algorithm nor did we see any benefits from using different iterative methods. Since the Hessian is (also for all practical purposes, in computation) a symmetric positive definite operator we prefer to use the PCG method over, e.g., GMRES. Having said this, we note that CLAIRE supports different Krylov subspace methods via \texttt{PETSc}~\cite{petsc-user-ref,petsc-web-page}. We discuss this in greater detail in the next subsection.

Since the considered optimization problem is, in general, non-convex, one additional challenge that arises is that the Hessian is not guaranteed to be positive definite, especially far away from a (local) minimizer. One approach to address this issue is to terminate the PCG algorithm as soon as one detects negative curvature. In this case, we use the former iterate of the PCG algorithm as a search direction. We consider a Gauss--Newton approximation to $\di{H}^{(k)}$~\cite{Mang:2015a,Mang:2017a,Mang:2019a} as an alternative to this approach. This approximation is guaranteed to be positive semi-definite. On the downside, we can expect the convergence to drop from quadratic to superlinear. This Gauss--Newton approximation is obtained by dropping all terms that involve the dual variable $\lambda$ in~\eqref{e:matvec} and~\eqref{e:incadj}, respectively. Notice that the final condition for the dual variable $\lambda$ in~\eqref{e:adjoint} corresponds to the mismatch between the transported intensities of the template image $m_0$ and the reference image $m_1$. Thus, as we approach a (local) minimizer of our problem, we can expect that $\lambda$ tends to zero; our Gauss--Newton approximation becomes exact and we recover quadratic convergence.

To further amortize computational costs, we do not invert $\di{H}^{(k)}$ exactly. We consider an inexact scheme~\cite{Dembo:1983a,Eisenstat:1996a,Nocedal:2006a}. This is accomplished by selecting the stopping condition for the PCG method to be proportional to the norm of the reduced gradient; as we approach a (local) minimizer, the tolerance decreases and we solve for the search direction more accurately. That is, we terminate the algorithm if
\[
\|\di{r}^{(k)} \|_\infty \leq \eta^{(k)} \|\di{g}^{(k)}\|_\infty, \quad \di{r}^{(k)} \defeq \di{H}^{(k)} \di{\tilde{v}}^{(k)}  + \di{g}^{(k)},
\]

\noindent with forcing sequence $\eta^{(k)} = \min (1/2, \sqrt{\|\di{g}^{(k)}\|_\infty})$ or $\eta^{(k)} = \min (1/2, \|\di{g}^{(k)}\|_\infty)$ for superlinear or quadratic convergence, respectively. See \algref{a:inner-iteration}, line \ref{alg:forcing} in the appendix.

In~\cite{Mang:2015a} we demonstrate that the preconditioned gradient descent scheme is less effective than our (Gauss--)Newton--Krylov scheme. As such, we only consider our \mbox{(Gauss--)}Newton--Krylov algorithm here.

We terminate the optimization if we reduce the gradient by $\epsilon_{\text{opt}} > 0$, i.e.,
\[
\|\di{g}^{(k)}\|_{\infty} \leq \epsilon_{\text{opt}}\|\di{g}^{(0)}\|_{\infty}
\]

\noindent or if $\|\di{g}^{(k)}\|_{\infty} \leq \scinum{1e-6}$. We have implemented alternative stopping criteria~\cite{Mang:2015a} but do not consider them here.

\subsubsection{Preconditioning}

The main cost of the {(Gauss--)}Newton--Krylov algorithm is the solution of the linear system
\begin{equation}\label{e:reducedspacekkt}
\di{H}^{(k)} \di{\tilde{v}}^{(k)} = -\di{g}^{(k)},\quad k = 1,2,3, \ldots,
\end{equation}

\noindent at each outer iteration $k$, with $\di{H}^{(k)} = \di{H}_{\text{reg}}^{(k)} + \di{H}_{\text{dat}}^{(k)}$, where $\di{H}_{\text{reg}}^{(k)} = \di{H}_{\text{reg}}\in\mathbb{R}^{dn,dn}$ is a discrete representation of the regularization operator $\op{H}_{\text{reg}} = \alpha\op{L}$ and $\di{H}_{\text{dat}}^{(k)} \in \mathbb{R}^{dn,dn}$ is the discrete version of $\op{H}_{\text{dat}}$ in \eqref{e:matvec}. For the model outlined in \secref{s:formulation} the Hessian behaves like a compact operator; large eigenvalues are associated with smooth eigenvectors and the eigenvectors become more oscillatory as the eigenvalues decrease~\cite{Mang:2015a}.

To amortize the computational costs of our algorithm and make it competitive with gradient descent schemes that consider first-order derivative information only, we have to design effective methods for preconditioning the linear system given above. That is, we seek a matrix $\di{M}^{(k)}\succ \di{0}$ such that, ideally, $(\di{M}^{(k)})^{-1} \di{H}^{(k)} \approx \di{I}_{dn}$, $\di{I}_{dn} \defeq \operatorname{diag}(1,\ldots,1) \in \mathbb{R}^{dn,dn}$. This makes approximations to $\di{H}^{(k)}$ (that are ``easy'' to invert) an obvious choice.

\paragraph{Regularization Preconditioner}

A common choice in PDE-constrained optimization is to consider the regularization operator $\di{H}_{\text{reg}}$ as a preconditioner $\di{M}$~\cite{Bui:2012a,BuiThanh:2013a,Alexanderian:2016a}. The preconditioned Hessian is a perturbation of the identity, i.e.,
\[
\big(\di{H}_{\text{reg}}\big)^{-1} \di{H}^{(k)}
= \di{I}_{dn} + \big(\di{H}_{\text{reg}}\big)^{-1}\di{H}_{\text{dat}}^{(k)}
\]

\noindent with $\di{I}_{dn} \defeq \operatorname{diag}(1,\ldots,1) \in \mathbb{R}^{dn,dn}$. Since $\di{H}_{\text{reg}} \succeq \di{0}$ is a (high-order) differential operator (typically, a Helmholtz type operator), its inverse acts as a smoother. We note that applying the inverse of $\di{H}_{\text{reg}}$ has a complexity of $\mathcal{O}(n \log n)$ in our implementation, i.e., we have to compute two FFTs and a diagonal scaling in the spectral domain using the appropriate Fourier coefficients. As such, this strategy for preconditioning the reduced space Hessian has vanishing costs. This preconditioner has been considered in~\cite{Mang:2015a,Mang:2016a,Mang:2016b,Mang:2018a}. The performance of this preconditioner is mesh independent (assuming we can entirely resolve the problem on the coarsest mesh). However, it deteriorates significantly as we decrease the regularization parameter $\alpha$.

\paragraph{Two-Level Preconditioner}

Inspired by multi-grid approaches, we designed a two-level preconditioner for the reduced space Hessian~\cite{Mang:2017a,Mang:2019a}. We use a coarse grid approximation of the inverse of the reduced space Hessian as a preconditioner. The basic idea is to iterate only on the low-frequency part and ignore the high-frequency components. That is, we use the inverse of the reduced space Hessian $\di{H}^{(k)}$, inverted on a coarser grid, as a preconditioner. This idea is motivated by the work in~\cite{Adavani:2008a,Biros:2008a,Giraud:2006a,Kaltenbacher:2003a,Kaltenbacher:2001a,King:1990a}. For simplicity of notation, we drop the dependence of the Hessian on the outer iteration index $k$.

We decompose the Hessian into two operators $\di{H}_L$ and $\di{H}_H$---one acting on low and the other acting on high frequencies, respectively. We denote the operators that project on the low and high-frequency subspaces by $\di{P}_L : \mathbb{R}^{dn} \to \mathbb{R}^{dn}$ and $\di{P}_H : \mathbb{R}^{dn} \to \mathbb{R}^{dn}$. Let $\di{e}_j\in \mathbb{R}^n$, $(\di{e}_j)_i = 1$ if $j=i$ and $(\di{e}_j)_i = 0$ for $i\not=j$, $i,j=1,\ldots,n$, denote an eigenvector of $\di{H}$ with $(\di{P}_L\di{H}\di{P}_H)\di{e}_j = (\di{P}_H\di{H}\di{P}_L)\di{e}_j = \di{0}$. Then, with $\di{P}_H + \di{P}_L = \di{I}_{dn}$, $ \di{I}_{dn} \defeq \operatorname{diag}(1,\ldots,1) \in \mathbb{R}^{dn,dn}$, we have
\[
\di{H}\di{e}_j
= (\di{P}_H + \di{P}_L)\di{H}(\di{P}_H + \di{P}_L)\di{e}_j
= \di{P}_H\di{H}\di{P}_H\di{e}_j + \di{P}_L\di{H}\di{P}_L\di{e}_j,
\]

\noindent In general, this equality will not hold. However, we are not interested in using this model as a surrogate for the Hessian $\di{H}$; we are merely interested in designing an effective preconditioner $\di{M}$ so that $\operatorname{cond}(\di{M}^{-1} \di{H}) \ll \operatorname{cond}(\di{H})$.

Suppose we can decompose $\di{\tilde{v}} \in \mathbb{R}^{dn}$ into a smooth component $\di{\tilde{v}}_L \in \mathbb{R}^{dn}$ and a high-frequency component $\di{\tilde{v}}_H \in \mathbb{R}^{dn}$, where each of these vectors can be found by solving
\[
\di{H}_L\di{\tilde{v}}_L = (\di{P}_L\di{H}\di{P}_L)\di{\tilde{v}}_L = -\di{P}_L\di{g}
\qquad\text{and}\qquad
\di{H}_H\di{\tilde{v}}_H = (\di{P}_H\di{H}\di{P}_H)\di{\tilde{v}}_H = -\di{P}_H\di{g},
\]

\noindent respectively. We use this construction to design an effective preconditioner for the smooth spectrum of our problem. Let $\di{r}\in\mathbb{R}^{dn}$ denote the vector we apply the inverse of our preconditioner $\di{M} \in \mathbb{R}^{dn,dn}$ to. Since our implementation is matrix-free, we iteratively solve $\di{M} \di{s} = \di{r}$ to obtain the action of the inverse of $\di{M}$ on $\di{r}$. In the spirit of the conceptual idea introduced above, we assume that we can decompose $\di{s}$ into a smooth component $\di{s}_L$ and a high-frequency component $\di{s}_H$. Let $\di{Q}_R \in \mathbb{R}^{dn/2,dn}$ denote a restriction operator and $\di{Q}_P \in \mathbb{R}^{dn,dn/2}$ denote prolongation operator. Moreover, let $\di{F}_L \in \mathbb{R}^{dn,dn}$ and  $\di{F}_H \in \mathbb{R}^{dn,dn}$ denote a low and high-pass filter, respectively. We project the vector $\di{r} \in \mathbb{R}^{dn}$ to a vector $\di{r}_L \in \mathbb{R}^{dn/2}$ by filtering the high-frequency components and restricting the resulting vector to a coarser mesh, i.e., $\di{r}_L =  \di{Q}_R \di{F}_L \di{r}$. Subsequently, we obtain the smooth component $\di{s}_L$ by solving
\[
\di{\tilde{M}}_L \di{\tilde{s}}_L = \di{Q}_R \di{F}_L \di{r}
\]

\noindent where $\di{\tilde{M}}_L \in \mathbb{R}^{dn/2,dn/2}$ is a coarse grid approximation of the low-frequency part of the reduced space Hessian $\di{H}$ and $\di{\tilde{s}}_L \in \mathbb{R}^{dn/2}$. This allows us to precondition the smooth part of $\di{r}$. We note that we do not precondition the high-frequency components $\di{r}_H \defeq \di{F}_H \di{r}$ of $\di{r}$, where $\di{F}_H \in \mathbb{R}^{dn,dn}$ is a high-pass filter with $\di{F}_H + \di{F}_L = \di{I}_{dn}$. Consequently, $\di{s}_H = \di{F}_H \di{r}$. In summary, the solution of $\di{M} \di{s} = \di{r}$ is given by
\[
\di{s} =
\di{s}_L + \di{s}_H
\approx
\di{Q}_P \di{F}_L \di{\tilde{s}}_L + \di{F}_H \di{r}
\approx
\di{Q}_P \di{F}_L (\di{\tilde{M}}_L)^{-1}\di{Q}_R \di{F}_L \di{r} + \di{F}_H \di{r}.
\]

\noindent To counter the fact that we leave the high-frequency components untouched, we do not directly apply this preconditioner to the reduced-space KKT system in \eqref{e:reducedspacekkt} but the regularization preconditioned system

\[
(\di{I} + \di{H}_{\text{reg}}^{\nicefrac{-1}{2}} \di{H}_{\text{data}}\di{H}_{\text{reg}}^{\nicefrac{-1}{2}}) \di{w} = - \di{H}_{\text{reg}}^{\nicefrac{-1}{2}} \di{g}
\]

\noindent where $\di{w} \defeq \di{H}_{\text{reg}}^{\nicefrac{1}{2}} \di{\tilde{v}}$. Notice that the square root of the inverse of $\di{H}_{\text{reg}}$ acts as a smoother. This scheme can be viewed as an approximation of a two-level multigrid V-cycle with an explicit (algebraic) smoother $\di{H}_{\text{reg}}^{\nicefrac{-1}{2}}$.

Before we explore extensions of this idea, we present some implementation aspects. We use spectral restriction and prolongation operators $\di{Q}_R$ and $\di{Q}_P$~\cite{Mang:2019a,Mang:2017a}. The operators $\di{F}_L$ and $\di{F}_H$ are implemented as cut-off filters in the frequency domain~\cite{Mang:2019a,Mang:2017a}. For the implementation of the coarse grid operator $\di{\tilde{M}}_L$ we have two choices. First, we can use a Galerkin discretization, which is formally given by $\di{\tilde{M}}_L = \di{Q}_R\di{H}\di{Q}_P$~\cite{Briggs:2000a}. The drawback of this approach is that we do not significantly reduce the computational costs compared to inverting the fine-grid Hessian, since each matvec necessitates the solution of the incremental state and adjoint equation at full resolution. Conversely, we can directly discretize the Hessian on a coarse grid to obtain $\di{\tilde{M}}_L$. This makes the implementation slightly more involved but reduces the computational costs drastically. We opt for the latter approach~\cite{Mang:2019a,Mang:2017a}.

To invert the matrix $\di{\tilde{M}}_L$ we have several options. Again, traditional direct methods are out of the question. However, we can use a nested Krylov-subspace method. If we use a Krylov-subspace method as an outer method (i.e., for computing the search direction), we have to select a tolerance for the inner Krylov-subspace method that is a fraction of the tolerance used to compute the search direction. Alternatively, we can replace the solver for the Newton step with a flexible Krylov-subspace method~\cite{Axelsson:1991a,Notay:2000a} and use a fixed number of iterations for the nested (inner) Krylov-subspace method. Alternatively, we can use a semi-iterative Chebyshev method~\cite{Gutknecht:2002a} with a fixed number of iterations on the inside. This yields a fixed linear operator for a particular choice of eigenvalue bounds~\cite{Golub:1961a}. These bounds can be estimated using a Lanczos method. We have tested and compared these approaches in~\cite{Mang:2017a,Mang:2019a}. This also includes the use of different Krylov-subspace methods for not only applying the preconditioner but also solving for the Newton step such as the standard and flexible GMRES method, the standard and flexible PCG method, or the Chebyshev method (some of which have been mentioned above). In~\cite{Mang:2019a}, we observed that the nested PCG method converges most quickly in the 3D setting.

\paragraph{Zero Velocity Approximation}

The preconditioner introduced in the former section requires a repeated evaluation of the incremental state and adjoint equations. The savings come from discretizing the reduced space Hessian on a mesh of half the resolution. In~\cite{Brunn:2020a} we developed a preconditioner that does not require solving any PDEs; the Hessian operator is fixed across all iterations. This is accomplished by fixing $\vec{v}$ to $\vec{v} = \vec{0}$ (our initial guess for the optimization problem). Under the assumption, the state equation simplifies to $\partial_t m = 0$, i.e., $m(t, \vec{x}) = m_0(\vec{x})$ for all $\vec{x} \in \Omega$ and $t\in[0,1]$. Likewise, we have $\partial_t \lambda = 0$, i.e., $\lambda(t, \vec{x}) = -(m_0(\vec{x}) - m_1(\vec{x}))$ for all $\vec{x} \in \Omega$ and $t\in[0,1]$. Inserting these expressions into the incremental state equation we have $\partial_t \tilde{m} = - \igrad m_0 \cdot \vec{\tilde{v}}$, which implies that $\tilde{m}(t=1) = - \igrad m_0 \cdot \vec{\tilde{v}}$. The incremental adjoint equation for the Gauss--Newton approximation for $\vec{v} = \vec{0}$ is given by $\partial_t \tilde{\lambda} = 0$, i.e., $\tilde{\lambda}(t, \vec{x}) = \igrad m_0(\vec{x}) \cdot \vec{\tilde{v}}(\vec{x})$ for all $\vec{x} \in \Omega$ and $t\in[0,1]$. Consequently, the Gauss--Newton approximation of the Hessian matvec for $\vec{v} = \vec{0}$ is given by
\[
\op{H}_{0}[\vec{\tilde{v}}](\vec{x}) = \alpha\op{L}\vec{\tilde{v}}(\vec{x}) + (\igrad m_0(\vec{x}) \otimes \igrad m_0(\vec{x})) \vec{\tilde{v}}(\vec{x}).
\]

This approximation deteriorates as we move away from our initial guess $\di{v}^{(0)} = \di{0}$. As a remedy, we replace $m_0$ in the expression above with our current estimate $m$ at $t=1$ at each outer iteration $k$ for a trial velocity $\di{v}^{(k)}$. Like in previous sections, we do not form or store $\di{H}_0 \in \mathbb{R}^{dn,dn}$ (the discrete version of $\op{H}_0$); we invert the matrix iteratively using a nested PCG method. To further reduce the computational costs, we combine the $\di{H}_0$ approximation with the two-level scheme discussed above. That is, we replace the coarse grid preconditioner $\di{\tilde{M}}_L$ with a coarse grid approximation of  $\op{H}_{0}$.

\subsection{Parameter Selection}\label{s:paratuning}

Based on empirical observations, we fix most of our numerical parameters. For the number of time steps $n_t$ in the numerical time integration we found that $n_t = 4$ provides sufficient accuracy to obtain a good matching between images at resolutions at the order of $256^3$ (standard size for brain images acquired in clinical practice). We set the tolerance for the relative reduction of the gradient (stopping condition for optimization) to $\epsilon_{\text{opt}} = \scinum{5e-2}$. We use a superlinear forcing sequence to compute the tolerance for the outer PCG algorithm. We use a two-level implementation of the zero velocity approximation of the reduced space Hessian as a preconditioner. The tolerance for the inner PCG to invert $\di{H}_0$ is 10 times smaller than the outer tolerance of the PCG. The formulation we consider for diffeomorphic image registration is an extension of what we discussed so far; it considers near-incompressible velocities. We describe this formulation is greater detail in the appendix. The regularization parameter for the $H^1$ penalty for the divergence of the velocity field is fixed and set to $\beta = \scinum{1e-4}$. We compute an optimal regularization parameter $\alpha$ as described below.

Several methods exist to estimate an optimal regularization parameter for inverse problems (see, e.g.,~\cite{Vogel:2002a} for examples). All of these methods have in common that the estimation of an optimal regularization parameter is expensive. Methods that assume that the differences between model output and observed data are random (such as, e.g., generalized cross validation) are not necessarily reliable in the context of image registration, since imaging noise is prone to be highly structured~\cite{Haber:2006a}. In our work, we consider a binary search for identifying an optimal value $\alpha$ for the regularization model for the velocity field $\vec{v}$~\cite{Mang:2015a,Mang:2019a,Himthani:2022a}. This approach is in spirit similar to an L-curve strategy. Related parameter continuation strategies have been considered in~\cite{Azencott:2010a,Haber:2000a,Haber:2006a}. As a measure for optimality, we select bounds on the determinant of the deformation gradient $\det \igrad \vec{y}$. Notice that we do not compute $\vec{y}$ to obtain this quantity but solve a transport problem (see appendix). Assuming that we start from an identity map $\operatorname{id}_{\mathbb{R}^d}$ the initial value for $\det \igrad \vec{y}$ is one (this is equivalent to $\vec{v} = \vec{0}$ in our formulation). Consequently, we assume that the map is diffeomorphic if $\det \igrad \vec{y} \geq 0$. This motivates the use of a lower bound $\epsilon_D > 0$. Since the determinant of the deformation gradient of $\vec{y}$ is inversely proportional to $\det \igrad \vec{y}^{-1}$, we use $1/\epsilon_D$ as an upper bound. Consequently, we require
\begin{equation}\label{e:ddgradbound}
\epsilon_D < \det \igrad \vec{y} < 1/\epsilon_D
\end{equation}

\noindent for any admissible $\vec{y}$. Our approach is as follows: We start with a regularization parameter of $\alpha = 1$ and reduce $\alpha$ by one order of magnitude until the condition in~\eqref{e:ddgradbound} is violated. Subsequently, we perform a binary search in the interval between the last value $\alpha$ for which~\eqref{e:ddgradbound} held and the value for which~\eqref{e:ddgradbound} was violated. For each new trial parameter $\alpha^{(l)}$ at level $l \in \mathbb{N}$, we use the control variable $\vec{v}_{\alpha}^{(l-1)} \defeq \vec{v}(\alpha^{(l-1)})$ obtained for $\alpha^{(l-1)}$ at the prior level $l-1$ as an initial guess to speed up convergence. More details can be found in~\cite{Mang:2015a,Mang:2019a}.

Obviously, this search is expensive since we have to solve the inverse problem for each trial $\alpha^{(l)}$, $l = 0,1,2,3, \ldots$. Once we have identified an adequate regularization parameter $\alpha^\star$ for a particular application, we perform a parameter continuation to speed up convergence. That is, we solve the inverse problem consecutively for different values for $\alpha$, starting with $\alpha^{(0)} = 1$ and subsequently reducing $\alpha^{(l)}$ by one order of magnitude until we reach the order of $\alpha^\star$. Then, we solve the problem one last time for $\alpha^\star$. For high regularization parameters $\alpha^{(l)}$ we essentially solve a convex problem; we expect quick convergence to a (local) minimizer. We use the estimate for the control variable $\vec{v}_{\alpha}^{(l-1)}$ as an initial guess for the next solve at level $l$. This does not significantly affect the runtime compared to directly solving our problem for $\alpha^\star$. Moreover, it ``convexifies'' the problem; we anticipate to more quickly converge to a ``better'' (local) minimizer and/or avoid getting trapped in ``less optimal'' local minima. We have compared this strategy against multi-scale (scale continuation) and multi-resolution (grid continuation) approaches in~\cite{Mang:2019a}. We observed the parameter continuation to be more stable and overall more effective. Combining parameter continuation with scale and/or grid continuation remains subject to future work.

Lastly, we note that machine learning has also recently been considered for regularization operator and parameter tuning~\cite{Wang:2020a,Al:2021a}.

\section{Results}\label{s:results}

We consider a slightly more involved formulation than the one presented in \secref{s:formulation}. We provide additional details in the appendix. We refer to~\cite{Mang:2016a,Gholami:2017a,Mang:2019a,Brunn:2020a} for weak and strong scaling results of our CPU and GPU implementation of CLAIRE, respectively. In the present work, we limit the performance analysis to a single GPU.

\subsection{Data}

\begin{figure}
\centering
\includegraphics[width=\textwidth]{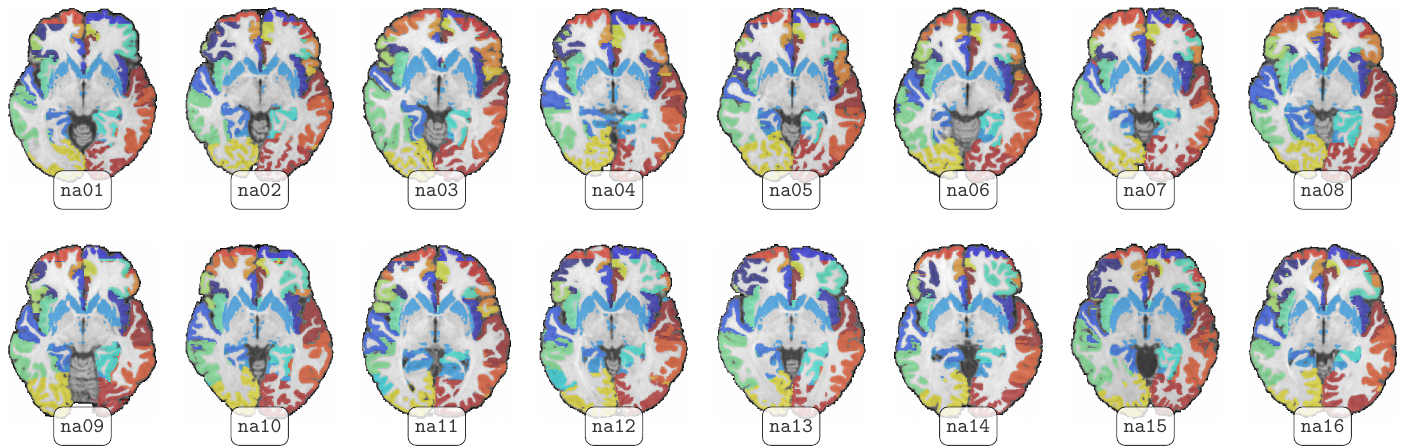}
\caption{{\tt NIREP} data repository~\cite{Christensen:2006a}. We show an axial slice of each dataset (slice number 128). The repository contains 16 rigidly aligned T1-weighted MRI brain datasets ({\tt na01}--{\tt na16})  of size $256\times300\times256$ voxels of different individuals. Each dataset is equipped with 32 labels of anatomical gray matter regions. We overlay these regions in different colors on the MRI data. We refer to~\cite{Christensen:2006a} for additional information about the datasets, the imaging protocol, and the preprocessing.}
\label{f:nirep-data}
\end{figure}

We report results for the {\tt NIREP} dataset~\cite{Christensen:2006a}. We refer to~\cite{Christensen:2006a} for additional information about the datasets, the imaging protocol, and the preprocessing. This repository contains 16 rigidly aligned T1-weighted MRI brain datasets ({\tt na01}--{\tt na16})  of size $256\times300\times256$ voxels of different individuals. Consequently, we invert for $3(256^2)(300) = \inum{58982400}$ unknowns. Each dataset is equipped with 33 labels for anatomical gray matter regions. These labels allow us to assess the performance of the registration; we assess registration accuracy by how well these labels are mapped to one another. To do so, we compute the so-called Dice between the label maps. For a Dice of one, the labels are in perfect agreement. For a Dice of zero, they do not overlap. Notice that the registration software does not consider these labels; registration is solely based on matching corresponding image intensities. That is, we do not explicitly minimize the alignment of the labels but only the mismatch between the data. We show the considered data in \figref{f:nirep-data}. In particular, we show axial slices of all 16 datasets with the associated labels in the overlay.

\subsection{Preconditioning}

\begin{figure}
\centering
\includegraphics[height=3cm]{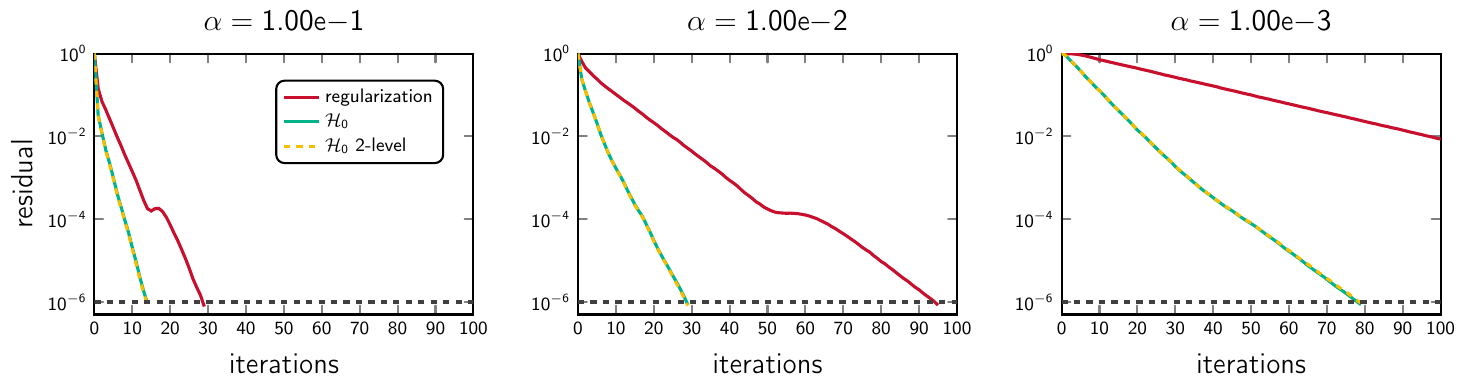}
\caption{Convergence of the PCG method for solving for the Newton step. We solve for the search direction at the solution of the registration problem (dataset {\tt na02} registered to {\tt na01}). We consider the squared $L^2$-distance. We solve this problem at the original resolution of the data. We consider three different preconditioners: The regularization preconditioner, the zero-velocity approximation of the reduced space Hessian, and the 2-level implementation of the zero-velocity approximation of the reduced space Hessian. We solve the problem for three different regularization parameter values (from left to right): $\alpha = \scinum{1e-1}$, $\alpha = \scinum{1e-2}$, and $\alpha = \scinum{1e-3}$. The tolerance for the PCG method is $\scinum{1e-6}$. We plot the relative residual.}
\label{f:convergence-preconditioners}
\end{figure}

We show representative results for the convergence of different preconditioners in \figref{f:convergence-preconditioners}. We consider the regularization preconditioner as well as two variants of the zero-velocity preconditioner---inverting the zero-velocity approximation of the reduced space Hessian on the fine mesh and a two-level implementation of this preconditioner. To test the performance, we invert the reduced space Hessian at the true solution of our problem. That is, we solve the registration problem between two images (dataset {\tt na02} registered to {\tt na01}) in our case. We then use the obtained velocity as iterate at which we compute the search direction. We set the tolerance for the PCG method to $\scinum{1e-6}$. We consider a squared $L^2$-distance as a similarity measure. We report results for the full resolution, only; $\vec{n}_x = (256,300,256)$. We report convergence results for three different choices of $\alpha$; $\alpha = \scinum{1e-1}$, $\alpha = \scinum{1e-2}$, and $\alpha = \scinum{1e-3}$, respectively.

The most important observations are: \begin{enumerate*}[label=(\roman*)] \item The convergence of all methods is sensitive with respect to the choice of $\alpha$. \item The zero-velocity approximation yields an improved rate of convergence. \item For the zero-velocity approximation, the convergence does not deteriorate as we switch from full resolution to a coarse resolution (2-level implementation)\end{enumerate*}.

\subsection{Regularization Parameter Search}

We set the regularization parameter for the divergence of the velocity to $\beta = \scinum{1e-5}$ and search for an optimal regularization parameter $\alpha$ using the scheme described in \secref{s:paratuning}. We register each image with all other images. We also perform the reverse registration. This results in a total of $16(15) = 240$ registrations. We consider a squared $L^2$-distance for the similarity measure.

\begin{figure}
\centering
\includegraphics[height=3cm]{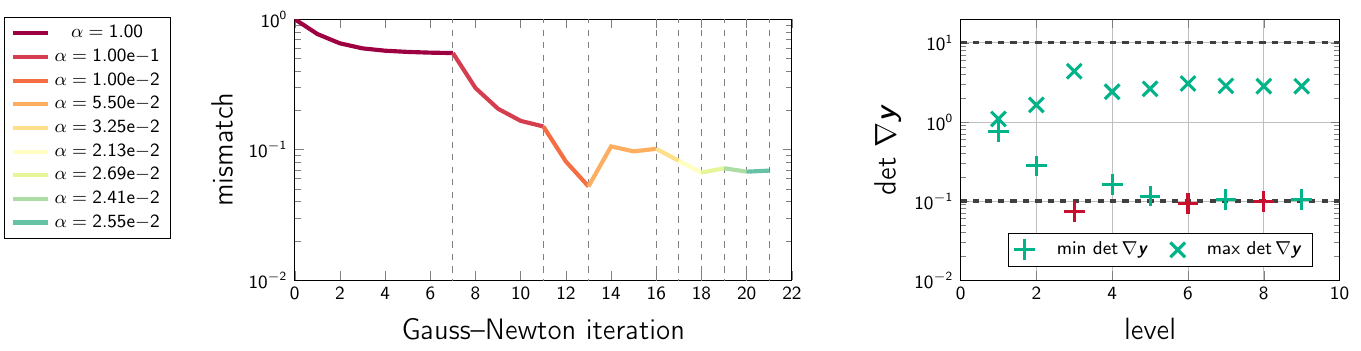}
\includegraphics[height=3cm]{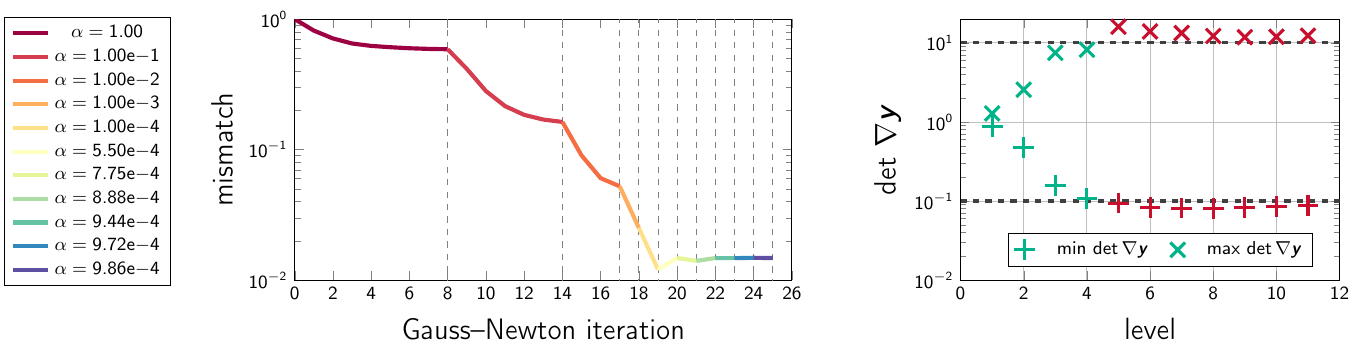}
\caption{Illustration of the parameter search for the registration of the dataset with id {\tt na14} to the dataset with id {\tt na01} (top panel) and in the reverse direction (bottom panel). The registration results are shown in \figref{f:claire-train-h1sdiv-betaw1e-5-bound-0d10} (bottom panel). We show (for each case) the trend of the mismatch for each choice of regularization parameter $\alpha$ (left) and the largest and smallest value of the determinant of the deformation gradient. In the latter plot we also show the lower and upper bound of 0.1 and 10, respectively, for the determinant of the deformation gradient as a dashed line. Whenever these bounds are violated, the marker switches from ``green'' to ``red.''. For the run shown in the top panel, the optimal regularization parameter is \snum{2.546875e-02}. For the run at the bottom, the optimal regularization parameter is $\scinum{1e-3}$.}
\label{f:claire-train-h1sdiv-betaw1e-5-bound-0d10-search}
\end{figure}

\begin{figure}
\centering
\includegraphics[width=\textwidth]{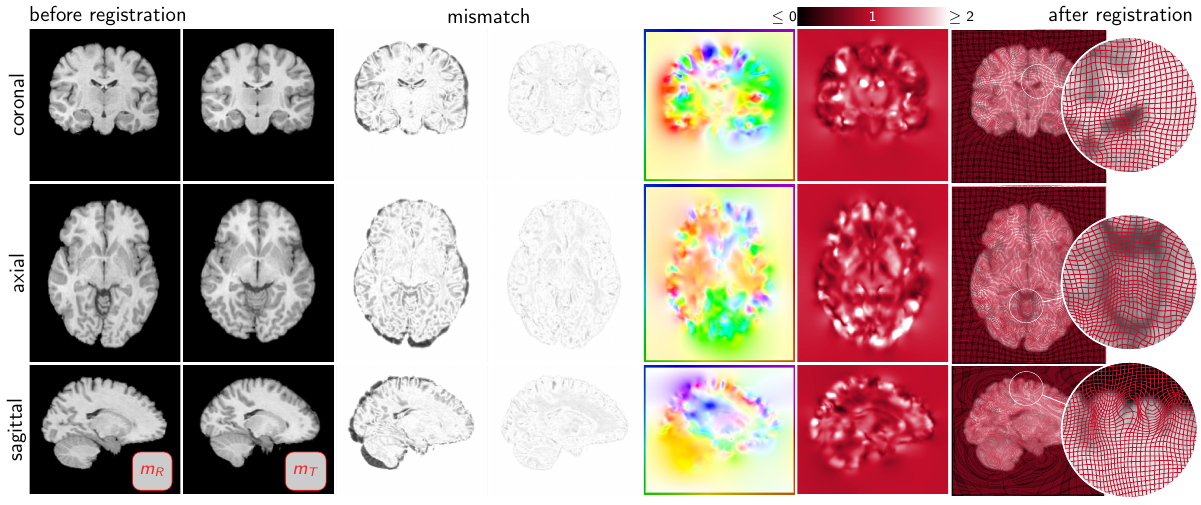}
\includegraphics[width=\textwidth]{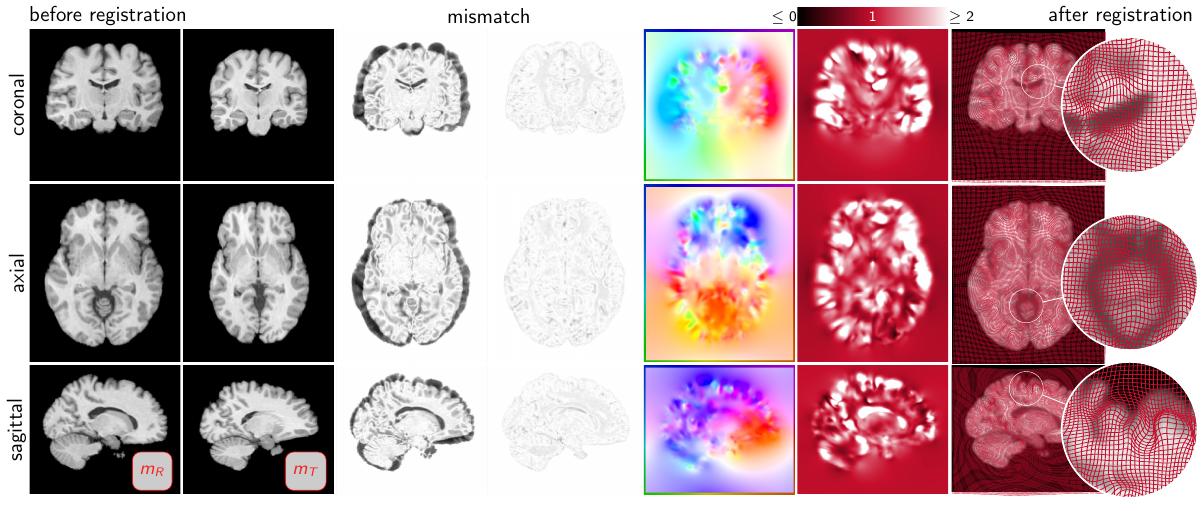}
\caption{Representative registration results for CLAIRE. We consider an $H^1$-seminorm as a regularization model for the  velocity field and an $H^1$-norm to regularize the divergence of the velocity field. We model near incompressible flows. The regularization parameter for the divergence is set to $\beta=\scinum{1e-5}$. The regularization for the velocity is estimated. The bound for the determinant of the deformation gradient is set to $\scinum{1e-1}$. We register the dataset with id {\tt na06} to the dataset with id {\tt na02} (top panel) and the dataset with id {\tt na14} to the dataset with id {\tt na01} (bottom panel) of the NIREP repository. The data is rigidly aligned. For each panel, we show the following: The top row shows the coronal view, the middle row the axial view, and the bottom row the sagittal view of the 3D volume. The columns are (from left to right) ($i$) the template image $m_0$, ($ii$) the reference image $m_1$, ($iii$) the residual differences between the reference image and the template image (before registration; large differences are colored in black and no residual difference are colored in white), ($iv$) the residual differences between the deformed template image and the reference image (after registration), ($v$) an illustration of the velocity field (color represents orientation; see boundary), ($vi$) visualization of the determinant of the deformation gradient (color bar on top), ($vii$) and an illustration of the projection of the computed deformation map onto the corresponding plane.\label{f:claire-train-h1sdiv-betaw1e-5-bound-0d10}}
\end{figure}

\begin{figure}
\centering
\includegraphics[height=3cm]{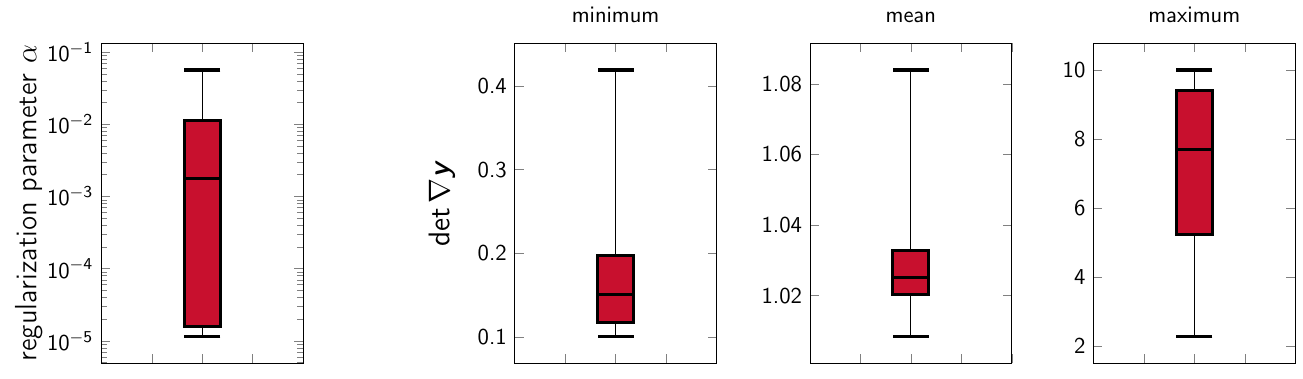}
\caption{Statistics for the estimation of the regularization parameter $\alpha$ across 240 registrations between all {\tt NIREP} datasets. The data has been rigidly registered. We report the estimated regularization parameters $\alpha$ (right) and the statistics for the minimum, mean and maximum of the determinant of the deformation gradient across each individual registration. The estimated regularization parameter $\alpha$ is $\snum{7.525404e-03}$ with with a standard deviation of $\snum{1.098077e-02}$, a median of $\snum{1.773437e-03}$, a minimal value of $\snum{1.140625e-05}$, and a maximal value of $\snum{5.640625e-02}$ across all 240 registrations. The 25th percentile is $\snum{1.562500e-05}$ and the 75th percentile is $\snum{1.140625e-02}$. The minimal, mean and maximal values are $\snum{1.667834e-01}$ (standard deviation: $\snum{6.118895e-02}$), $\snum{1.027673e+00}$ (standard deviation: \snum{1.182927e-02}), and  $\snum{7.149775e+00}$ (standard deviation: \snum{2.323663e+00}).}
\label{f:claire-train-h1sdiv-betaw1e-5-bound-0d10-estimated-parameters}
\end{figure}

We illustrate the search for an optimal regularization parameter for two registration problems ({\tt na01} to {\tt na14} and {\tt na14} to {\tt na01}) in \figref{f:claire-train-h1sdiv-betaw1e-5-bound-0d10-search}. We show representative registration results for two images from the considered {\tt NIREP} dataset in \figref{f:claire-train-h1sdiv-betaw1e-5-bound-0d10}. We report statistics for the estimated regularization parameter $\alpha$ in \figref{f:claire-train-h1sdiv-betaw1e-5-bound-0d10-estimated-parameters} (left plot). We also compute the minimal, mean, and maximum value of the determinant of the deformation gradient for all registrations. We report the statistics across all 240 registrations for these in \figref{f:claire-train-h1sdiv-betaw1e-5-bound-0d10-estimated-parameters} (plots to the right). For the minimum value of the determinant of the deformation gradient, we obtained $\snum{1.667834e-01}$ with with a standard deviation of $\snum{6.118895e-02}$, an overall lowest minimum value of $\snum{1.000982e-01}$ and an overall largest minimum value of $\snum{4.192143e-01}$. For the mean value of the determinant of the deformation gradient, we obtained $\snum{1.027673e+00}$ with with a standard deviation of $\snum{1.182927e-02}$, an overall lowest minimum value of $\snum{1.008390e+00}$ and an overall largest minimum value of $\snum{1.083912e+00}$. For the maximum value of the determinant of the deformation gradient, we obtained $\snum{7.149775e+00}$ with a standard deviation of \snum{2.323663e+00}, an overall lowest minimum value of $\snum{2.268732e+00}$ and an overall largest minimum value of $\snum{9.993294e+00}$. We report the workload for this search in \tabref{t:workload-parametersearch}.

\begin{table}
\caption{Workload for estimating the regularization parameter $\alpha$. We consider a squared $L^2$-distance measure. All reported numbers are computed across all 240 registrations. We report the number of outer iterations, the number of Hessian matvecs, the number of PDE solves, the relative mismatch after registration, the relative change of the norm of the gradient, and the runtime (in seconds). These numbers are for solving the inverse problem multiple times; we search for an optimal regularization parameter using a binary search.}
\label{t:workload-parametersearch}
\centering
\begin{tabular}{lrrrrrrr}
\toprule
% & \quad\bf mean & \quad\bf stdev & \quad\bf min & \quad\bf max & \quad\bf median & \quad\bf 1st QT & \quad\bf 3rd QT \\
& \bf mean & \bf stdev & \bf min & \bf max & \bf median & \bf 1st QT & \bf 3rd QT \\
\midrule
iterations &            $\fnum{2.053750e+01}$ & \fnum{2.268266e+00} & \inum{1.700000e+01} & \inum{2.800000e+01} & \inum{2.000000e+01} & \inum{1.900000e+01} & \inum{2.200000e+01} \\
matvecs    &            $\fnum{1.669708e+02}$ & \fnum{1.247491e+02} & \inum{3.900000e+01} & \inum{3.880000e+02} & \inum{9.150000e+01} & \inum{6.400000e+01} & \inum{3.260000e+02} \\
PDE solves &            $\fnum{4.479750e+02}$ & \fnum{3.188766e+02} & \inum{1.340000e+02} & \inum{1.008000e+03} & \inum{2.520000e+02} & \inum{1.900000e+02} & \inum{8.665000e+02} \\
mismatch   & \color{red}$\snum{4.508358e-02}$ & \snum{3.992633e-02} & \snum{4.421833e-03} & \snum{1.756399e-01} & \snum{3.182786e-02} & \snum{1.000444e-02} & \snum{7.095301e-02} \\
gradient   &            $\snum{1.393934e-02}$ & \snum{5.732892e-03} & \snum{3.166478e-03} & \snum{3.350037e-02} & \snum{1.338130e-02} & \snum{8.995437e-03} & \snum{1.835874e-02} \\
runtime    & \color{red}$\fnum{1.513511e+01}$ & \fnum{1.132311e+01} & \fnum{3.974425e+00} & \fnum{3.786656e+01} & \fnum{8.124641e+00} & \fnum{5.876944e+00} & \fnum{2.923212e+01} \\
\bottomrule
\end{tabular}
\end{table}

The most important observations are: \begin{enumerate*}[label=(\roman*)] \item We can efficiently determine an adequate regularization parameter with an average runtime of $\fnum{1.513511e+01}$ seconds (standard deviation: \fnum{1.132311e+01} seconds), \item the computed deformation maps are diffeomorphic (up to numerical accuracy), and \item we overall obtain high-quality registration results with precise control on the determinant of the deformation gradient\end{enumerate*}.

\subsection{Registration Accuracy}

In this section we assess the registration accuracy. In particular, we report the Dice values for the parameter search described in the former section. Aside from considering a squared $L^2$-distance we also report registration accuracy for normalized cross correlation as a similarity measure (see appendix for details). In \figref{f:dice-individual-labels} we report the Dice score for the individual labels. We report the statistics for the 240 registration runs in \tabref{t:dice-union-stats}. Here, we compute the union of all 33 labels and report the global Dice score.  We report additional results in the appendix.

\begin{figure}
\centering
\includegraphics[width=\textwidth]{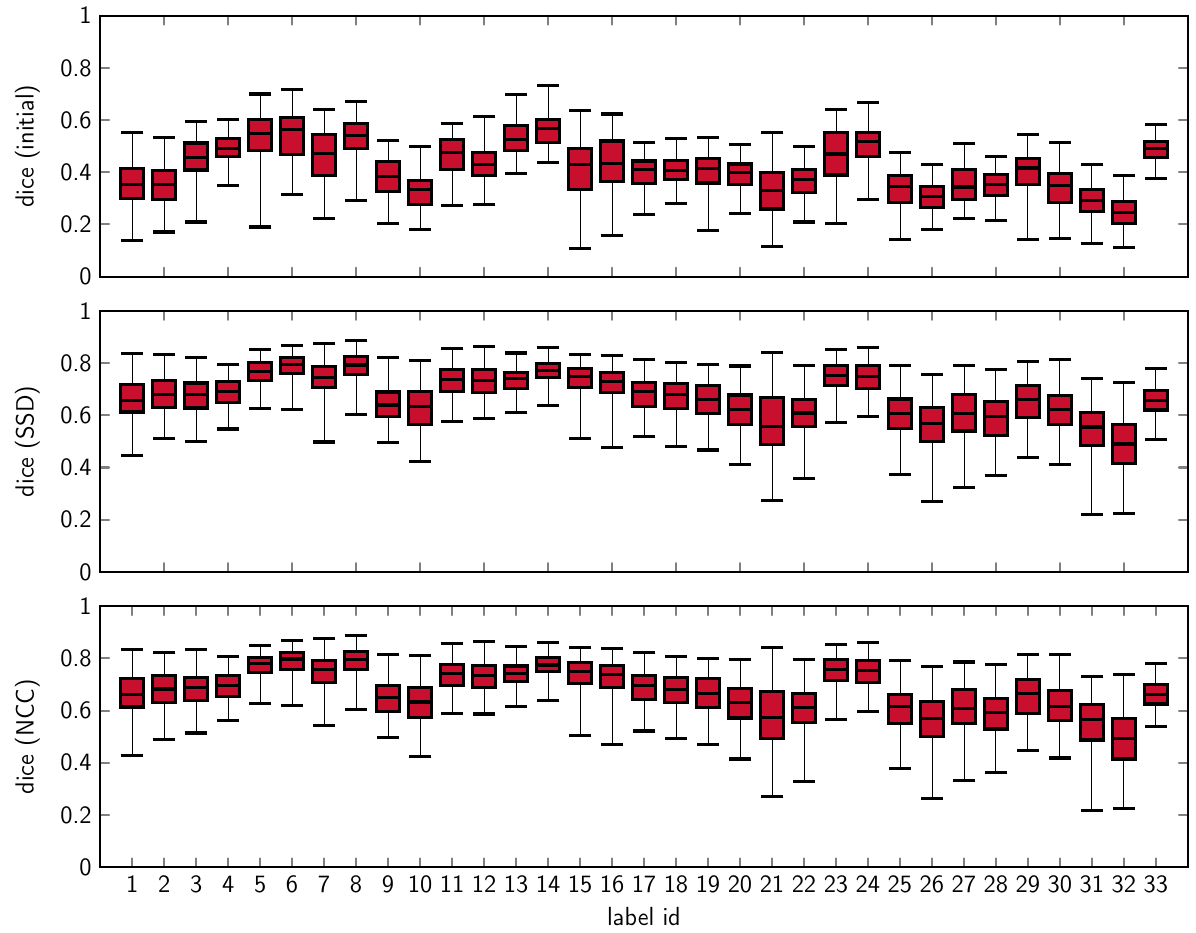}
\caption{Average Dice score for individual labels. We show box plots for the Dice score for each individual label. The statistics are computed for all 240 registrations. The top row corresponds to the Dice values before registration. The middle row shows values for the Dice score after registration using a squared $L^2$-distance as a similarity measure. The bottom row shows the results obtained for normalized cross correlation. We report statistics for the union of these labels in \tabref{t:dice-union-stats}.}
\label{f:dice-individual-labels}
\end{figure}

The most important observations are: \begin{enumerate*}[label=(\roman*)] \item CLAIRE yields an excellent agreement for the overall Dice with an increase from $\fnum{5.506221e-01}$ (standard deviation: \fnum{4.128437e-02})  before registration to $\fnum{8.307170e-01}$ (standard deviation: \fnum{5.912268e-02}) for the squared $L^2$-distance and $\fnum{8.346387e-01}$ (standard deviation: \fnum{5.996745e-02}) for normalized cross correlation. \item The performance for the squared $L^2$-distance and normalized cross correlation are en par for our current implementation\end{enumerate*}.

\begin{table}
\caption{Average DICE values. We report the mean, min, max, and median value as well as the 1st quantile and the 3rd quantile. These values are computed for the union of all labels. We report the initial values in the first row. The values after diffeomorphic registration based on the squared $L^2$-distance and normalized cross correlation are reported in the second and third rows, respectively. We report the scores for the individual labels in \figref{f:dice-individual-labels}.}
\label{t:dice-union-stats}
\centering
\begin{tabular}{lrrrrrrr}
\toprule
%& \quad\bf mean & \quad\bf stdev & \quad\bf min & \quad\bf max & \quad\bf median & \quad\bf 1st quantile & \quad\bf 3rd quantile \\
& \bf mean & \bf stdev & \bf min & \bf max & \bf median & \bf 1st quantile & \bf 3rd quantile \\
\midrule
initial & $\fnum{5.506221e-01}$ & \fnum{4.128437e-02} & \fnum{4.209550e-01} & \fnum{6.246068e-01} & \fnum{5.546036e-01} & \fnum{5.270924e-01} & \fnum{5.825034e-01} \\
SSD     & $\fnum{8.307170e-01}$ & \fnum{5.912268e-02} & \fnum{6.966414e-01} & \fnum{9.221054e-01} & \fnum{8.419006e-01} & \fnum{7.850278e-01} & \fnum{8.841150e-01} \\
NCC     & $\fnum{8.346387e-01}$ & \fnum{5.996745e-02} & \fnum{6.994966e-01} & \fnum{9.229853e-01} & \fnum{8.437425e-01} & \fnum{7.843089e-01} & \fnum{8.894927e-01} \\
\bottomrule
\end{tabular}
\end{table}

\subsection{Convergence and Runtime}

\begin{figure}
\centering
\includegraphics[width=0.4\textwidth]{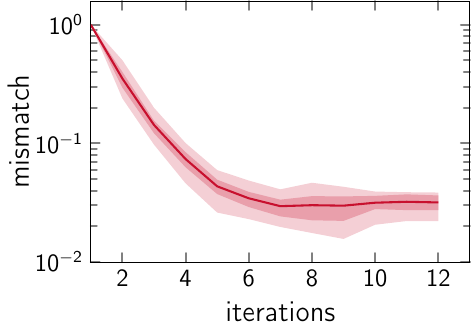}
\caption{Convergence behavior. We plot the relative reduction of the mismatch versus the number of iterations. The plot is generated for 240 registrations. We execute the algorithm for a regularization parameter value of $\alpha = \snum{1.773437e-03}$. We solve the problem without using any continuation scheme. The average runtime is \fnum{3.615090e+00}. The solid line represents the mean convergence for the data mismatch. We also show the envelopes for the 25th to 75th quantile and the 5th to the 95th quantile for the values of the mismatch. The average number of iterations is $\fnum{1.055417e+01}$. We show the trend until iteration 12.}
\label{f:claire-exec-nirep-ssd-convergence}
\end{figure}

In the former section, we have seen how CLAIRE performs when searching for an optimal regularization parameter for each individual volume. In the current section, we fix the regularization parameter to the mean optimal value of $\alpha = \snum{1.773437e-03}$ determined in the former section and focus on computational performance. We plot the residual vs. the number of outer iterations in \figref{f:claire-exec-nirep-ssd-convergence}. Here, we solve the inverse problem for a fixed $\alpha = \snum{1.773437e-03}$ without performing any scale, grid, or parameter continuation. We report the runtime for our parameter continuation scheme for a target regularization parameter $\alpha = \snum{1.773437e-03}$ in \tabref{t:workload-parametercontinuation}.

\begin{table}
\caption{Workload for estimating the regularization parameter $\alpha$. We consider a squared $L^2$-distance measure. All reported numbers are computed across all 240 registration. We report the number of outer iterations, the number of Hessian matvecs, the number of PDE solves, the relative mismatch after registration, the relative change of the norm of the gradient, and the runtime (in seconds). These numbers are for solving the inverse problem multiple times; we search for an optimal regularization parameter using a binary search.}
\label{t:workload-parametercontinuation}
\centering
\begin{tabular}{lrrrrrrr}
\toprule
%& \quad\bf mean & \quad\bf stdev & \quad\bf min & \quad\bf max & \quad\bf median & \quad\bf 1st QT & \quad\bf 3rd QT \\
& \bf mean & \bf stdev & \bf min & \bf max & \bf median & \bf 1st QT & \bf 3rd QT \\
\midrule
iterations &            $\fnum{1.055417e+01}$ & \fnum{1.246410e+00} & \inum{9.000000e+00} & \inum{1.500000e+01} & \inum{1.000000e+01}  & \inum{1.000000e+01} & \inum{1.100000e+01} \\
matvecs    &            $\fnum{2.452500e+01}$ & \fnum{5.659017e+00} & \inum{1.400000e+01} & \inum{4.700000e+01} & \inum{2.400000e+01}  & \inum{2.100000e+01} & \inum{2.750000e+01} \\
PDE solves &            $\fnum{8.215833e+01}$ & \fnum{1.366393e+01} & \inum{5.800000e+01} & \inum{1.360000e+02} & \inum{8.200000e+01}  & \inum{7.400000e+01} & \inum{8.900000e+01} \\
mismatch   & \color{red}$\snum{4.398711e-02}$ & \snum{1.036629e-02} & \snum{1.905666e-02} & \snum{8.086897e-02} & \snum{4.292485e-02}  & \snum{3.690673e-02} & \snum{5.008053e-02} \\
gradient   &            $\snum{3.940170e-02}$ & \snum{6.484699e-03} & \snum{2.134815e-02} & \snum{4.973159e-02} & \snum{4.012440e-02}  & \snum{3.623087e-02} & \snum{4.409452e-02} \\
runtime    & \color{red}$\fnum{3.838822e+00}$ & \fnum{6.482507e-01} & \fnum{2.703537e+00} & \fnum{6.398900e+00} & \fnum{3.779393e+00}  & \fnum{3.450854e+00} & \fnum{4.159088e+00} \\
\bottomrule
\end{tabular}
\end{table}

The most important observations are: \begin{enumerate*}[label=(\roman*)] \item We can solve the inverse problem in under 4 seconds (on average, the runtime is $\fnum{3.838822e+00}$ seconds; standard deviation: \fnum{6.482507e-01} seconds), with a minimum runtime of under 3 seconds and a maximum runtime of slightly above 6 seconds. \item We converge in about 12 iterations to a stable solution of our problem (the mismatch stagnates), where a majority of the runs we have executed converge after only $\fnum{1.055417e+01}$ (standard deviation: \fnum{1.246410e+00}). \item Once we have determined an adequate regularization parameter for a particular application, we can solve the problem quickly with an accuracy that is equivalent to the more expensive parameter search considered in the section above as judged by the relative reduction of the mismatch\end{enumerate*}.

\section{Conclusions}\label{s:conclusions}

We have reviewed our past work on scalable algorithms for diffeomorphic image registration. Several issues remain.

Our implementation currently only supports the registration of images acquired with the same imaging modality. Developing an effective solver for other distance measures remains subject to future work. We have worked on several numerical schemes for preconditioning the reduced space Hessian. The spectral preconditioner is extremely efficient to apply but its performance deteriorates as we reduce the regularization parameter. This is true for all other schemes we have implemented to precondition the reduced Hessian. Although they are more effective than the simple spectral preconditioner, developing a scheme that has a rate of convergence that is mesh-independent and at the same time independent of the choice of the regularization (parameter) remains subject to future work.

Our 3D GPU implementation currently only supports stationary velocities. These velocities do not define a proper metric in the Riemannian space of diffeomorphic flows. While we have implemented a MATLAB prototype version of a solver that supports time varying velocities, this implementation has not yet been ported to the C++ implementation of CLAIRE.

Another challenge in diffeomorphic image registration is how to handle data that underwent topological changes (e.g., the emergence of a tumor or tissue being removed due to clinical intervention). One possibility to handle this is to introduce additional biophysical constraints~\cite{Gholami:2017a,Gholami:2016a,Scheufele:2020a,Scheufele:2019a,Mang:2017c,Gooya:2012a,Mang:2020a,Hogea:2008a,Zacharaki:2008b,Zacharaki:2009a,Hogea:2008b}. On the downside, this makes the problem much more challenging to solve since we not only invert for a deformation map but also for the parameters of the model. More generic approaches to deal with changes in topology are described in~\cite{Hsieh:2022b,Li:2012a,Franccois:2022a,Antonsanti:2021a,Sukurdeep:2022a}.

\noindent{\bf Acknowledgements.} This work was in part supported by the National Science Foundation (NSF) through the grants DMS-2012825 and DMS-2145845. Any opinions, findings, and conclusions or recommendations expressed herein are those of the author and do not necessarily reflect the views of the NSF. This work was completed in part with resources provided by the Research Computing Data Core at the University of Houston. The author would like to thank George Biros, Malte Brunn, Amir Gholami, Naveen Himthani, Jae Youn Kim, and Miriam Schulte for their numerous contributions to this work.

\section*{Appendix}

\subsection*{PDE Constraints}

Below, we will revisit some of the problem formulations we have considered in our past work. These are extensions to the formulation considered in \secref{s:formulation}. We limit the description of our methodology to the most basic formulation for simplicity. The default formulation implemented in our current GPU version in CLAIRE is different~\cite{Brunn:2020a,Brunn:2021a}.

\subsubsection*{Non-Stationary Velocities}

In~\cite{Mang:2015a}, we consider stationary and non-stationary velocities. For non-stationary velocities, the reduced gradient in~\eqref{e:rgrad} is given by
\[
\vec{g}(\vec{v}) \defeq
\alpha \op{L} \vec{v} + \lambda \igrad m.
\]

\subsubsection*{(Near-)Incompressible Diffeomorphisms}

In~\cite{Mang:2015a}, we augment the formulation in~\eqref{e:varopt-claire} by introducing the incompressibility constraint $\idiv \vec{v} = 0$. A similar formulation has been considered in~\cite{Chen:2012a}. For the primal-dual optimal variables $(m^\star, \vec{v}^\star, \lambda^\star, \rho^\star)$, the associated KKT conditions are given by
\begin{subequations}\label{e:kkt-incompressible}
\begin{align}
\partial_t m^\star + \vec{v}^\star \cdot \igrad m^\star & = 0
&& \text{in}\;\; (0,1] \times \Omega,
\\
m^\star &= m_0
&& \text{in}\;\; \{0\} \times \Omega,
\\
-\partial_t \lambda^\star + \idiv \lambda^\star\vec{v}^\star  & = 0
&& \text{in}\;\; [0,1) \times \Omega,
\\
\lambda^\star &=  -(m^\star - m_1)
&& \text{in}\,\,\{1\} \times \Omega,
\\
\idiv \vec{v}^\star &= 0
&& \text{in}\,\,\Omega,\label{e:incompressible}
\\
\alpha \op{L} \vec{v}^\star + \igrad  \rho^\star + \int_0^1\lambda^\star \igrad m^\star\, \d t &= \vec{0}
&& \text{in}\,\,\Omega.
\end{align}
\end{subequations}

We eliminate the incompressibility constraint~\eqref{e:incompressible} and the dual variable $\rho$ from the optimality system stated above to obtain the expression
\[
\alpha \op{L} \vec{v}
+ \int_0^1\lambda \igrad m \d t
- \igrad \ilap^{-1} \idiv \int_0^1\lambda \igrad m\, \d t
\]

\noindent for the reduced gradient. The remaining PDE operators in~\eqref{e:kkt-incompressible} for $m$ and $\lambda$ in the associated KKT system are identical.

In~\cite{Mang:2016a}, we relaxed the incompressiblity constraint by introducing an additional control variable $w$ to obtain $\idiv \vec{v} = w$. This allows us to model near-incompressible deformations. After eliminating the constraint $\idiv \vec{v} = w$ and the associated dual variable $\rho$ from the KKT system, we obtain the reduced gradient
\[
\alpha \op{L} \vec{v}^\star
+ \int_0^1\lambda^\star \igrad m^\star \d t
- \igrad (\alpha(\beta (-\ilap^{-1} + \operatorname{id}))^{-1} + \operatorname{id})^{-1} \ilap^{-1} \idiv \int_0^1\lambda^\star \igrad m^\star\, \d t,
\]

\noindent Here, $\beta > 0$ denotes the regularization parameter of the regularizer for the second control variable $w$. We consider an $H^1$-norm. We refer to~\cite{Mang:2016a} for additional details. This represents the default model implemented in the hardware-accelerated implementation of CLAIRE~\cite{Mang:2019a,Brunn:2020a,Brunn:2021a}. The results reported in this study also consider this formulation. The regularization model for the velocity field is an $H^1$-seminorm.

Aside from this, we have also explored a model of incompressible flows that promotes shear~\cite{Mang:2016a}. To do so, we introduce a nonlinear regularization model. In particular, we replaced the regularization model for $\vec{v}$ by
\[
|\vec{v}|_{H^1(\Omega)}^{(1+\nu)/2\nu} =
\frac{2 \nu}{\nu + 1}
\int_{\Omega}
\left(
\mathcal{E}[\vec{v}] \colon \mathcal{E}[\vec{v}]
\right)^{(1 + \nu)/2\nu}\, \d \vec{x},
\]

\noindent where
\[
\mathcal{E}[\vec{v}] \defeq \frac{1}{2} \left((\igrad_d \vec{v}) + (\igrad_d \vec{v})^\mathsf{T}\right),
\quad\igrad_d \vec{v} \defeq
\left( \begin{array}{c} (\igrad v_1)^\mathsf{T} \\ \vdots \\(\igrad v_d)^\mathsf{T} \end{array}\right)
\in \mathbb{R}^{d,d},
\]

\noindent denotes the strain tensor, and $\nu > 0$ controls the non-linearity. With this regularization model in conjunction with the incompressibility constraint $\idiv \vec{v} = 0$ we obtain a Stokes-like optimality system with a viscosity that depends on the strain rate. The reduced gradient is given by
\[
- \operatorname{div}\left( 2 \operatorname{tr}(\mathcal{E}[\vec{v}]\mathcal{E}[\vec{v}] )^{(1-\nu)/2\nu} \mathcal{E}[\vec{v}]\right) + \igrad \rho + \int_0^1\lambda \igrad m \d t,
\]

\noindent where $\operatorname{div}( \vec{A} ) = \left(\idiv \vec{a}_1, \ldots, \idiv \vec{a}_d\right) \in \mathbb{R}^d$ for an arbitrary $d \times d$ matrix $\vec{A}$ with columns $\vec{a}_i \in \mathbb{R}^d$, $i = 1, \ldots, d$.  In the limit $\nu \to \infty$ this model behaves like total variation regularization. For $\nu \in (0,1)$ we obtain a shear thickening and for $\nu > 1$ a shear thinning fluid. Likewise to the linear case, we can eliminate the incompressiblity constraint and the associated dual variable $\rho$ from the optimality system. We refer to \cite{Mang:2016a} for additional details.

\subsubsection*{Optimal Transport}

In our past work, we have not only introduced new hard or soft constraints for $\vec{v}$ but also considered a different forward model for transporting $m$. In particular, we use the continuity equation
\[
\begin{aligned}
\partial_t m + \idiv m \vec{v} & = 0
&& \text{in}\;\; [0,1) \times \Omega,
\end{aligned}
\]

\noindent with initial condition $m = m_0$ in $\{0\} \times \Omega$ to model the transport of the intensities of the template image $m_0$. In this model, mass is conserved. This establishes a connection to optimal transport~\cite{Angenent:2003a,Rehman:2009a,Chen:2018b}. We refer to~\cite{Mang:2017a} for additional details.

\subsection*{Deformation Gradient}

In the context of image registration, the determinant of the deformation gradient $\det \igrad \vec{y}$ is often used to assess invertibility of $\vec{y}$ as well as a measure of local volume change in the context of morphometry and shape analysis. In the framework of continuum mechanics, we can obtain this information from the deformation tensor field $\vec{f} : [0,1] \times \bar{\Omega} \to \mathbb{R}^{d, d}$, where $\vec{f}$ is related to $\vec{v}$ by
\begin{equation}\label{e:def-grad}
\partial_t \vec{f} + (\vec{v} \cdot \igrad_d) \vec{f}
= (\igrad_d \vec{v}) \vec{f} \;\; \text{in} \;\; \Omega \times (0,1],
\qquad \vec{f} = \vec{I}_{d} \;\; \text{in} \quad \Omega \times \{0\},
\end{equation}

\noindent with periodic boundary conditions on $\partial\Omega$. Here, $\vec{I}_d = \operatorname{diag}(1,\ldots,1)\in\mathbb{R}^{d,d}$. In our implementation we use $\det \vec{f}_1$ with $\vec{f}_1 \defeq \vec{f}(t=1,\,\cdot\,)$ as a surrogate for $\det \igrad \vec{y}$,

\subsection*{Normalized Cross Correlation}

Aside from using the squared $L^2$-distance, we also consider normalized cross correlation as a distance measure. We note that we have not presented results for normalized cross correlation elswhere. The choice of the similarity measure in general only affects the final condition of the dual variable. The normalized cross correlation distance measure is given by
\begin{equation}
\operatorname{dist}_{\text{NCC}} (m(1),m_1)
= 1 - \frac{\langle m(1), m_1 \rangle^2_{L^2(\Omega)}}
{\langle m_1, m_1 \rangle_{L^2(\Omega)}
\langle m(1), m(1) \rangle_{L^2(\Omega)}},
\end{equation}

\noindent where
\[
\langle u, w \rangle_{L^2(\Omega)} = \int_{\Omega} u(\vec{x}) w(\vec{x})\, \d \vec{x}
\]

\noindent denotes the standard $L^2$-inner product on $\Omega \subset \mathbb{R}^d$ for arbitrary functions $u : \bar{\Omega} \to \mathbb{R}$, $w : \bar{\Omega} \to \mathbb{R}$. Using this distance, the final condition for the adjoint equation is given by
\[
\lambda(1,\vec{x})
= -2\frac{\langle m_1,m(1)\rangle_{L^2(\Omega)}}{\|m(1)\|^2_{L^2(\Omega)} \|m_R\|^2_{L^2(\Omega)}}
\left(\frac{\langle m_1,m(1)\rangle_{L^2(\Omega)}}{\|m(1)\|^2_{L^2(\Omega)}}m(1,\vec{x}) - m_1(\vec{x})\right).
\]

Similarly, the expression for the final condition of the incremental dual variable $\tilde{\lambda}$ is given by
\[
\tilde{\lambda}(1,\vec{x})
= \frac{2(q_1 m_1(\vec{x}) + q_2 m(1,\vec{x}) - q_3 \tilde{m}(1,\vec{x}))}{\|m_1\|^2_{L^2(\Omega)}},
\]

\noindent where
\[
\begin{aligned}
q_1 &=
2\frac{\langle m_1, m(1)\rangle_{L^2(\Omega)}
\langle m(1), \tilde{m}(1)\rangle_{L^2(\Omega)}}{\|m(1)\|^4_{L^2(\Omega)}}
- \frac{\langle m_R, \tilde{m}(1)\rangle_{L^2(\Omega)}}{\|m(1)\|^2_{L^2(\Omega)}},
\\
q_2 & =
4\frac{\langle m_R,m(1)\rangle ^2_{L^2(\Omega)} \langle m(1),\tilde{m}(1) \rangle_{L^2(\Omega)}}{\|m(1)\|^6_{L^2(\Omega)}}
-2\frac{\langle m_R, m(1)\rangle_{L^2(\Omega)} \langle m_R, \tilde{m}(1)\rangle_{L^2(\Omega)}}{\|m(1)\|^4_{L^2(\Omega)}},
\\
q_3 &= \frac{\langle m_R, m(1)\rangle^2_{L^2(\Omega)}}{\|m(1)\|^4_{L^2(\Omega)}}.
\end{aligned}
\]

\subsection*{Newton--Krylov Algorithm}

We summarize our Newton--Krylov algorithm here. The outer iterations are given in \algref{a:outer-iteration}. The inner iterations (i.e., the computation of the search direction) is given in \algref{a:inner-iteration}. We describe this algorithm in some detail in \secref{s:optimization}.

\begin{algorithm}
\caption{Inexact Newton--Krylov method (outer iterations). We use the relative norm of the reduced gradient with tolerance $\epsilon_{\text{opt}} > 0$ as stopping criterion.}
\label{a:outer-iteration}
\begin{algorithmic}[1]
\STATE{$k\leftarrow0$}
\STATE{initial guess $\di{v}^{(k)} \leftarrow \di{0}$}
\STATE{$\di{m}^{(k)} \leftarrow$ solve state equation in~\eqref{e:varopt:constraint} forward in time given $\di{v}^{(k)}$}
\STATE{$j^{(k)} \leftarrow$ evaluate objective functional~\eqref{e:varopt:objective} given $\di{m}^{(k)}$ and $\di{v}^{(k)}$}
\STATE{$\di{\lambda}^{(k)} \leftarrow$ solve adjoint equation~\eqref{e:adjoint} backward in time given $\di{v}^{(k)}$ and $\di{m}^{(k)}$}
\STATE{$\di{g}^{(k)} \leftarrow$ evaluate reduced gradient~\eqref{e:rgrad} given $\di{m}^{(k)}$, $\di{\lambda}^{(k)}$ and $\di{v}^{(k)}$}
\WHILE{$\|\di{g}^{(k)}\|_\infty > \|\di{g}^{(0)}\|_\infty\epsilon_{\text{opt}}$\label{alg:stop-newton}}
    \STATE{$\di{\tilde{v}}^{(k)}\leftarrow$ solve $\di{H}^{(k)} \di{\tilde{v}}^{(k)} = -\di{g}^{(k)}$ given $\di{m}^{(k)}$, $\di{\lambda}^{(k)}$, $\di{v}^{(k)}$, and $\di{g}^{(k)}$ (see \algref{a:inner-iteration})}
    \STATE{$\gamma^{(k)} \leftarrow$ perform line search on $\di{\tilde{v}}^{(k)}$ subject to Armijo condition}
    \STATE{$\di{v}^{(k+1)} \leftarrow \di{v}^{(k)} + \gamma^{(k)}\di{\tilde{v}}^{(k)}$}
    \STATE{$\di{m}^{(k+1)} \leftarrow$ solve state equation~\eqref{e:varopt:constraint} forward in time given $\di{v}^{(k+1)}$}
    \STATE{$j^{(k+1)} \leftarrow$ evaluate~\eqref{e:varopt:objective} given $\di{m}^{(k+1)}$ and $\di{v}^{(k+1)}$}
    \STATE{$\di{\lambda}^{(k+1)} \leftarrow$ solve adjoint equation~\eqref{e:adjoint} backward in time given $\di{v}^{(k+1)}$ and $\di{m}^{(k+1)}$}
    \STATE{$\di{g}^{(k+1)} \leftarrow$ evaluate~\eqref{e:rgrad} given $\di{m}^{(k+1)}$, $\di{\lambda}^{(k+1)}$ and $\di{v}^{(k+1)}$}
    \STATE{$k \leftarrow k + 1$}
\ENDWHILE
\end{algorithmic}
\end{algorithm}

\begin{algorithm}
\caption{Newton step (inner iterations). We illustrate the solution of the reduced KKT system~\eqref{e:reducedspacekkt} using a PCG method at a given outer iteration $k\in\mathbb{N}$. We use a superlinear forcing sequence to compute the tolerance $\eta^{(k)}$ for the PCG method (inexact solve).}
\label{a:inner-iteration}
\begin{algorithmic}[1]
\STATE{\textbf{input:} $\di{m}^{(k)}$, $\di{\lambda}^{(k)}$, $\di{v}^{(k)}$, $\di{g}^{(k)}$, $\di{g}^{(0)}$}
\STATE{$l \leftarrow 0$}
\STATE{set $\epsilon_H \leftarrow\min\big(0.5,\|\di{g}^{(k)}\|_\infty^{1/2}\big)$, \quad $\di{\tilde{v}}^{(l)} \leftarrow \di{0}$, \quad $\di{r}^{(l)} \leftarrow - \di{g}^{(k)}$ \label{alg:init}}
\STATE{$\di{z}^{(l)} \leftarrow$ apply preconditioner $\di{M}^{-1}$ to $\di{r}^{(l)}$}
\STATE{$\di{s}^{(l)} \leftarrow \di{z}^{(l)}$}
\WHILE{$l < n$}
    \STATE{$\di{\tilde{m}}^{(l)} \leftarrow$ solve~\eqref{e:incstate} forward in time given $\di{m}^{(k)}$, $\di{v}^{(k)}$ and $\di{\tilde{v}}^{(l)}$\label{alg:incstate}}
    \STATE{$\di{\tilde{\lambda}}^{(l)} \leftarrow$ solve~\eqref{e:incadj} backward in time given $\di{\lambda}^{(k)}$, $\di{v}^{(k)}$, $\di{\tilde{m}}^{(l)}$ and $\di{\tilde{v}}^{(l)}$\label{alg:inc-adjoint}}
    \STATE{$\di{\tilde{s}}^{(l)} \leftarrow$ apply $\di{H}^{(l)}$ to $\di{s}^{(l)}$ given $\di{\lambda}^{(k)}$, $\di{m}^{(k)}$, $\di{\tilde{m}}^{(l)}$ and $\di{\tilde{\lambda}}^{(l)}$ (see~\eqref{e:matvec})}
    \STATE{$\kappa^{(l)} \leftarrow\langle\di{r}^{(l)},\di{z}^{(l)}\rangle/\langle\di{s}^{(l)},\di{\tilde{s}}^{(l)}\rangle$}
    \STATE{$\di{\tilde{v}}^{(l+1)} \leftarrow \di{\tilde{v}}^{(l)} + \kappa^{(l)}\di{s}^{(l)}$}
    \STATE{$\di{r}^{(l+1)}\leftarrow\di{r}^{(l)}-\kappa^{(l)}\di{\tilde{s}}^{(l)}$}
    \STATE{\textbf{if} $\|\di{r}^{(l+1)}\|_2 < \epsilon_H$ \textbf{break}\label{alg:forcing}}
    \STATE{$\di{z}^{(l+1)} \leftarrow$ apply preconditioner $\di{M}^{-1}$ to $\di{r}^{(l+1)}$}
    \STATE{$\mu^{(l)} \leftarrow \langle\di{z}^{(l+1)},\di{r}^{(l+1)}\rangle/\langle\di{z}^{(l)},\di{r}^{(l)}\rangle$}
    \STATE{$\di{s}^{(l+1)} \leftarrow \di{z}^{(l+1)} + \mu^{(l)}\di{s}^{(l)}$}
    \STATE{$l \leftarrow l + 1$}
\ENDWHILE
\STATE{\textbf{output:} $\di{\tilde{v}}^{(k)} \leftarrow\di{\tilde{v}}^{(l+1)}$}
\end{algorithmic}
\end{algorithm}

\subsection*{Hardware}

We execture CLAIRE on the Sabine Cluster of the Research Computing Data Core at the University of Houston. Sabine hosts a total of 5704\,CPU cores in 169\,compute and 12\,GPU nodes. We limit the experiments to our GPU implementation. The associated nodes are equipped with a Intel Xeon E5-2680v4 CPUs (2 sockets with 28 cores) with 256\,GB of memory. Each node is also equiped with 8 NVIDIA V100 GPUs with a total of 40,960 cores and 128\,GB of memory.

\subsection*{Additional Results}

We report more detailed results for the registration accuracy of CLAIRE in this section. The statistics for the Dice for the squared $L^2$-distance with respect to each individual label is reported in \tabref{t:dice-individual-labels-l2}. The associated results for normalized cross correlation are reported in \tabref{t:dice-individual-labels-ncc}. These results are for the parameter search for 240 registration (all-to-all) of the NIREP dataset.

\begin{table}
\caption{Average DICE values. We report the mean, min, max, and median value as well as the 1st quantile and the 3rd quantile. These values are computed for each individual label for all 240 registrations. These results are obtained for the squared $L^2$-distance measure.\label{t:dice-individual-labels-l2}}
\centering
\begin{tabular}{rrrrrrrr}
\toprule
%\bf label id & \quad\bf mean & \quad\bf stdev & \quad\bf min & \quad\bf max & \quad\bf median & \quad\bf 1st quantile & \quad\bf 3rd quantile \\
\bf label id & \bf mean & \bf stdev & \bf min & \bf max & \bf median & \bf 1st quantile & \bf 3rd quantile \\
\midrule
  1 & $\fnum{6.614078e-01}$ & \fnum{8.163262e-02} & \fnum{4.449623e-01} & \fnum{8.380036e-01} & \fnum{6.577743e-01}  & \fnum{6.125849e-01} & \fnum{7.192794e-01} \\
  2 & $\fnum{6.777377e-01}$ & \fnum{7.035507e-02} & \fnum{5.101232e-01} & \fnum{8.334061e-01} & \fnum{6.806626e-01}  & \fnum{6.290541e-01} & \fnum{7.322086e-01} \\
  3 & $\fnum{6.777758e-01}$ & \fnum{6.725452e-02} & \fnum{4.996823e-01} & \fnum{8.223014e-01} & \fnum{6.814115e-01}  & \fnum{6.282392e-01} & \fnum{7.240040e-01} \\
  4 & $\fnum{6.887056e-01}$ & \fnum{5.406500e-02} & \fnum{5.475523e-01} & \fnum{7.940540e-01} & \fnum{6.911166e-01}  & \fnum{6.501439e-01} & \fnum{7.307965e-01} \\
  5 & $\fnum{7.650083e-01}$ & \fnum{4.588336e-02} & \fnum{6.262010e-01} & \fnum{8.520660e-01} & \fnum{7.691985e-01}  & \fnum{7.332846e-01} & \fnum{8.025228e-01} \\
  6 & $\fnum{7.887633e-01}$ & \fnum{4.185679e-02} & \fnum{6.238973e-01} & \fnum{8.674311e-01} & \fnum{7.964776e-01}  & \fnum{7.591309e-01} & \fnum{8.224897e-01} \\
  7 & $\fnum{7.384833e-01}$ & \fnum{6.988593e-02} & \fnum{4.976362e-01} & \fnum{8.759184e-01} & \fnum{7.439456e-01}  & \fnum{7.065561e-01} & \fnum{7.872277e-01} \\
  8 & $\fnum{7.860385e-01}$ & \fnum{5.388558e-02} & \fnum{6.018148e-01} & \fnum{8.876697e-01} & \fnum{7.928674e-01}  & \fnum{7.560349e-01} & \fnum{8.259297e-01} \\
  9 & $\fnum{6.418489e-01}$ & \fnum{6.607327e-02} & \fnum{4.963126e-01} & \fnum{8.200293e-01} & \fnum{6.394636e-01}  & \fnum{5.959484e-01} & \fnum{6.896039e-01} \\
 10 & $\fnum{6.294914e-01}$ & \fnum{8.114641e-02} & \fnum{4.227619e-01} & \fnum{8.114497e-01} & \fnum{6.335764e-01}  & \fnum{5.648403e-01} & \fnum{6.909873e-01} \\
 11 & $\fnum{7.322581e-01}$ & \fnum{5.527108e-02} & \fnum{5.772586e-01} & \fnum{8.575132e-01} & \fnum{7.364543e-01}  & \fnum{6.913351e-01} & \fnum{7.741084e-01} \\
 12 & $\fnum{7.289749e-01}$ & \fnum{5.821253e-02} & \fnum{5.881687e-01} & \fnum{8.650798e-01} & \fnum{7.329974e-01}  & \fnum{6.877279e-01} & \fnum{7.744374e-01} \\
 13 & $\fnum{7.339755e-01}$ & \fnum{4.576029e-02} & \fnum{6.106043e-01} & \fnum{8.388674e-01} & \fnum{7.417891e-01}  & \fnum{7.020168e-01} & \fnum{7.632993e-01} \\
 14 & $\fnum{7.683593e-01}$ & \fnum{4.228743e-02} & \fnum{6.374640e-01} & \fnum{8.611016e-01} & \fnum{7.709381e-01}  & \fnum{7.444253e-01} & \fnum{7.980782e-01} \\
 15 & $\fnum{7.362278e-01}$ & \fnum{5.672607e-02} & \fnum{5.100584e-01} & \fnum{8.345714e-01} & \fnum{7.480395e-01}  & \fnum{7.064810e-01} & \fnum{7.788159e-01} \\
 16 & $\fnum{7.212034e-01}$ & \fnum{6.146189e-02} & \fnum{4.753525e-01} & \fnum{8.277491e-01} & \fnum{7.307394e-01}  & \fnum{6.874794e-01} & \fnum{7.625749e-01} \\
 17 & $\fnum{6.798459e-01}$ & \fnum{6.218946e-02} & \fnum{5.198067e-01} & \fnum{8.128076e-01} & \fnum{6.898382e-01}  & \fnum{6.336941e-01} & \fnum{7.271112e-01} \\
 18 & $\fnum{6.730072e-01}$ & \fnum{6.885251e-02} & \fnum{4.797340e-01} & \fnum{8.039358e-01} & \fnum{6.805469e-01}  & \fnum{6.246541e-01} & \fnum{7.208454e-01} \\
 19 & $\fnum{6.550582e-01}$ & \fnum{7.280819e-02} & \fnum{4.669929e-01} & \fnum{7.943333e-01} & \fnum{6.596829e-01}  & \fnum{6.067236e-01} & \fnum{7.150986e-01} \\
 20 & $\fnum{6.207931e-01}$ & \fnum{7.920572e-02} & \fnum{4.119435e-01} & \fnum{7.891933e-01} & \fnum{6.228308e-01}  & \fnum{5.659685e-01} & \fnum{6.778063e-01} \\
 21 & $\fnum{5.694698e-01}$ & \fnum{1.219740e-01} & \fnum{2.723240e-01} & \fnum{8.410739e-01} & \fnum{5.564289e-01}  & \fnum{4.875263e-01} & \fnum{6.690745e-01} \\
 22 & $\fnum{6.018655e-01}$ & \fnum{8.139789e-02} & \fnum{3.587843e-01} & \fnum{7.901492e-01} & \fnum{6.089306e-01}  & \fnum{5.576285e-01} & \fnum{6.598770e-01} \\
 23 & $\fnum{7.495623e-01}$ & \fnum{5.436848e-02} & \fnum{5.713219e-01} & \fnum{8.513940e-01} & \fnum{7.541657e-01}  & \fnum{7.146508e-01} & \fnum{7.924217e-01} \\
 24 & $\fnum{7.462623e-01}$ & \fnum{5.346887e-02} & \fnum{5.961795e-01} & \fnum{8.613484e-01} & \fnum{7.472308e-01}  & \fnum{7.026412e-01} & \fnum{7.920293e-01} \\
 25 & $\fnum{6.044068e-01}$ & \fnum{8.180108e-02} & \fnum{3.749306e-01} & \fnum{7.916583e-01} & \fnum{6.078357e-01}  & \fnum{5.489038e-01} & \fnum{6.626674e-01} \\
 26 & $\fnum{5.617581e-01}$ & \fnum{9.126121e-02} & \fnum{2.705264e-01} & \fnum{7.580618e-01} & \fnum{5.685580e-01}  & \fnum{5.012706e-01} & \fnum{6.289436e-01} \\
 27 & $\fnum{6.000844e-01}$ & \fnum{9.789865e-02} & \fnum{3.242557e-01} & \fnum{7.923154e-01} & \fnum{6.072877e-01}  & \fnum{5.399182e-01} & \fnum{6.806048e-01} \\
 28 & $\fnum{5.869787e-01}$ & \fnum{8.381110e-02} & \fnum{3.703028e-01} & \fnum{7.766931e-01} & \fnum{5.938315e-01}  & \fnum{5.243618e-01} & \fnum{6.519905e-01} \\
 29 & $\fnum{6.501103e-01}$ & \fnum{7.761405e-02} & \fnum{4.398960e-01} & \fnum{8.048787e-01} & \fnum{6.596296e-01}  & \fnum{5.922963e-01} & \fnum{7.143204e-01} \\
 30 & $\fnum{6.175218e-01}$ & \fnum{8.311267e-02} & \fnum{4.115361e-01} & \fnum{8.141801e-01} & \fnum{6.208778e-01}  & \fnum{5.636647e-01} & \fnum{6.754548e-01} \\
 31 & $\fnum{5.471829e-01}$ & \fnum{9.302607e-02} & \fnum{2.199474e-01} & \fnum{7.405272e-01} & \fnum{5.553624e-01}  & \fnum{4.845604e-01} & \fnum{6.113987e-01} \\
 32 & $\fnum{4.855354e-01}$ & \fnum{1.026684e-01} & \fnum{2.241389e-01} & \fnum{7.259933e-01} & \fnum{4.901026e-01}  & \fnum{4.153047e-01} & \fnum{5.664702e-01} \\
 33 & $\fnum{6.557103e-01}$ & \fnum{5.652823e-02} & \fnum{5.072552e-01} & \fnum{7.801358e-01} & \fnum{6.572667e-01}  & \fnum{6.202701e-01} & \fnum{6.965002e-01} \\
\bottomrule
\end{tabular}
\end{table}

\begin{table}
\caption{Average DICE values. We report the mean, min, max, and median value as well as the 1st quantile and the 3rd quantile. These values are computed for each individual label for all 240 registrations. The results are obtained for the normalized cross correlation distance measure.\label{t:dice-individual-labels-ncc}}
\centering
\small
\begin{tabular}{rrrrrrrr}
\toprule
%\bf label id & \quad\bf mean &\quad\bf stdev  & \quad\bf min & \quad\bf max & \quad\bf median & \bf 1st quantile & \bf 3rd quantile \\
\bf label id & \bf mean & \bf stdev  & \bf min & \bf max & \bf median & \bf 1st quantile & \bf 3rd quantile \\
\midrule
  1 & $\fnum{6.646245e-01}$ & \fnum{8.176533e-02} & \fnum{4.271562e-01} & \fnum{8.351873e-01} & \fnum{6.603834e-01}  & \fnum{6.139712e-01} & \fnum{7.236381e-01} \\
  2 & $\fnum{6.809139e-01}$ & \fnum{7.042589e-02} & \fnum{4.908037e-01} & \fnum{8.222979e-01} & \fnum{6.827985e-01}  & \fnum{6.315021e-01} & \fnum{7.348062e-01} \\
  3 & $\fnum{6.832128e-01}$ & \fnum{6.577472e-02} & \fnum{5.140685e-01} & \fnum{8.339254e-01} & \fnum{6.896432e-01}  & \fnum{6.387323e-01} & \fnum{7.274741e-01} \\
  4 & $\fnum{6.938833e-01}$ & \fnum{5.208379e-02} & \fnum{5.622606e-01} & \fnum{8.082243e-01} & \fnum{6.967105e-01}  & \fnum{6.545833e-01} & \fnum{7.343959e-01} \\
  5 & $\fnum{7.723829e-01}$ & \fnum{4.313417e-02} & \fnum{6.282176e-01} & \fnum{8.502115e-01} & \fnum{7.801983e-01}  & \fnum{7.466076e-01} & \fnum{8.038595e-01} \\
  6 & $\fnum{7.901338e-01}$ & \fnum{4.244970e-02} & \fnum{6.201017e-01} & \fnum{8.672997e-01} & \fnum{7.979042e-01}  & \fnum{7.589037e-01} & \fnum{8.236082e-01} \\
  7 & $\fnum{7.434875e-01}$ & \fnum{6.679160e-02} & \fnum{5.430817e-01} & \fnum{8.770417e-01} & \fnum{7.579916e-01}  & \fnum{7.066143e-01} & \fnum{7.914173e-01} \\
  8 & $\fnum{7.873854e-01}$ & \fnum{5.267087e-02} & \fnum{6.049122e-01} & \fnum{8.876202e-01} & \fnum{7.964604e-01}  & \fnum{7.559844e-01} & \fnum{8.272419e-01} \\
  9 & $\fnum{6.475919e-01}$ & \fnum{6.626134e-02} & \fnum{4.985244e-01} & \fnum{8.163170e-01} & \fnum{6.494139e-01}  & \fnum{5.974171e-01} & \fnum{6.955769e-01} \\
 10 & $\fnum{6.324533e-01}$ & \fnum{8.004244e-02} & \fnum{4.251186e-01} & \fnum{8.096926e-01} & \fnum{6.328206e-01}  & \fnum{5.734933e-01} & \fnum{6.880340e-01} \\
 11 & $\fnum{7.364940e-01}$ & \fnum{5.459386e-02} & \fnum{5.876998e-01} & \fnum{8.579182e-01} & \fnum{7.434320e-01}  & \fnum{6.959451e-01} & \fnum{7.783319e-01} \\
 12 & $\fnum{7.305973e-01}$ & \fnum{5.882237e-02} & \fnum{5.871212e-01} & \fnum{8.638319e-01} & \fnum{7.343781e-01}  & \fnum{6.899488e-01} & \fnum{7.728674e-01} \\
 13 & $\fnum{7.400326e-01}$ & \fnum{4.384538e-02} & \fnum{6.145643e-01} & \fnum{8.472025e-01} & \fnum{7.434464e-01}  & \fnum{7.098510e-01} & \fnum{7.709003e-01} \\
 14 & $\fnum{7.725713e-01}$ & \fnum{4.181219e-02} & \fnum{6.390662e-01} & \fnum{8.598188e-01} & \fnum{7.746713e-01}  & \fnum{7.494126e-01} & \fnum{8.027118e-01} \\
 15 & $\fnum{7.398256e-01}$ & \fnum{5.853863e-02} & \fnum{5.036516e-01} & \fnum{8.416118e-01} & \fnum{7.510435e-01}  & \fnum{7.028113e-01} & \fnum{7.856740e-01} \\
 16 & $\fnum{7.253655e-01}$ & \fnum{6.341164e-02} & \fnum{4.694308e-01} & \fnum{8.373636e-01} & \fnum{7.391979e-01}  & \fnum{6.892219e-01} & \fnum{7.716689e-01} \\
 17 & $\fnum{6.860176e-01}$ & \fnum{6.296884e-02} & \fnum{5.218480e-01} & \fnum{8.224599e-01} & \fnum{6.957854e-01}  & \fnum{6.415348e-01} & \fnum{7.351976e-01} \\
 18 & $\fnum{6.784766e-01}$ & \fnum{6.793262e-02} & \fnum{4.928209e-01} & \fnum{8.058791e-01} & \fnum{6.811836e-01}  & \fnum{6.315514e-01} & \fnum{7.272490e-01} \\
 19 & $\fnum{6.599681e-01}$ & \fnum{7.317841e-02} & \fnum{4.706312e-01} & \fnum{7.993696e-01} & \fnum{6.654087e-01}  & \fnum{6.106841e-01} & \fnum{7.219690e-01} \\
 20 & $\fnum{6.254474e-01}$ & \fnum{8.007284e-02} & \fnum{4.143991e-01} & \fnum{7.952639e-01} & \fnum{6.310412e-01}  & \fnum{5.717089e-01} & \fnum{6.843186e-01} \\
 21 & $\fnum{5.758028e-01}$ & \fnum{1.198489e-01} & \fnum{2.707386e-01} & \fnum{8.414808e-01} & \fnum{5.740407e-01}  & \fnum{4.946511e-01} & \fnum{6.724949e-01} \\
 22 & $\fnum{6.042834e-01}$ & \fnum{8.298360e-02} & \fnum{3.283801e-01} & \fnum{7.949468e-01} & \fnum{6.127728e-01}  & \fnum{5.531416e-01} & \fnum{6.668831e-01} \\
 23 & $\fnum{7.536164e-01}$ & \fnum{5.598206e-02} & \fnum{5.661990e-01} & \fnum{8.525112e-01} & \fnum{7.584206e-01}  & \fnum{7.157118e-01} & \fnum{7.957802e-01} \\
 24 & $\fnum{7.486471e-01}$ & \fnum{5.423909e-02} & \fnum{5.971644e-01} & \fnum{8.605013e-01} & \fnum{7.524305e-01}  & \fnum{7.072420e-01} & \fnum{7.905340e-01} \\
 25 & $\fnum{6.093552e-01}$ & \fnum{8.079757e-02} & \fnum{3.766774e-01} & \fnum{7.936598e-01} & \fnum{6.172543e-01}  & \fnum{5.497407e-01} & \fnum{6.606353e-01} \\
 26 & $\fnum{5.653397e-01}$ & \fnum{9.157560e-02} & \fnum{2.629715e-01} & \fnum{7.676776e-01} & \fnum{5.700334e-01}  & \fnum{5.013281e-01} & \fnum{6.359944e-01} \\
 27 & $\fnum{6.034608e-01}$ & \fnum{9.466128e-02} & \fnum{3.313133e-01} & \fnum{7.864196e-01} & \fnum{6.096703e-01}  & \fnum{5.498806e-01} & \fnum{6.816588e-01} \\
 28 & $\fnum{5.887141e-01}$ & \fnum{8.278953e-02} & \fnum{3.641458e-01} & \fnum{7.781930e-01} & \fnum{5.927602e-01}  & \fnum{5.297028e-01} & \fnum{6.464821e-01} \\
 29 & $\fnum{6.534107e-01}$ & \fnum{7.976854e-02} & \fnum{4.476718e-01} & \fnum{8.167070e-01} & \fnum{6.667361e-01}  & \fnum{5.890384e-01} & \fnum{7.192366e-01} \\
 30 & $\fnum{6.179502e-01}$ & \fnum{8.301268e-02} & \fnum{4.181951e-01} & \fnum{8.154287e-01} & \fnum{6.169072e-01}  & \fnum{5.626400e-01} & \fnum{6.784739e-01} \\
 31 & $\fnum{5.527666e-01}$ & \fnum{9.308986e-02} & \fnum{2.185331e-01} & \fnum{7.306268e-01} & \fnum{5.658186e-01}  & \fnum{4.875091e-01} & \fnum{6.219828e-01} \\
 32 & $\fnum{4.883446e-01}$ & \fnum{1.024081e-01} & \fnum{2.235305e-01} & \fnum{7.397722e-01} & \fnum{4.918621e-01}  & \fnum{4.146773e-01} & \fnum{5.698949e-01} \\
 33 & $\fnum{6.623968e-01}$ & \fnum{5.328180e-02} & \fnum{5.382469e-01} & \fnum{7.801456e-01} & \fnum{6.622045e-01}  & \fnum{6.251797e-01} & \fnum{6.992224e-01} \\
\bottomrule
\end{tabular}
\end{table}

\printbibliography

\end{document}